# To Get Overall Shapes and New Data of the 120-Cell and the 600-Cell

*For commemorating great mathematician H. S. M. Coxeter*

**Kaida Shi**

State Key Laboratory of CAD&CG, Zhejiang University,
Hangzhou 310027, Zhejiang Province, China
Department of Mathematics, Zhejiang Ocean University,
Zhoushan 316004, Zhejiang Province, China

**Abstract** This research will be helpful for people to display the 2-dimensional projective models of 4-variable actual problems in many fields, in order to investigate deeply those actual problems. By using the theory of *N*-dimensional finite rotation group of the regular polytopes, the author established the 2-dimensional projective model of 4-dimensional rectangular coordinate system, and deduced a transformation matrix, and adopt it to display successfully the 2-dimensional overall shapes of two most complicated regular polytopes 120-Cell and 600-Cell. In the meantime, the author calculated all the vertex coordinates and determine the joint relationships between adjacent vertices of the regular polytopes 120-Cell and 600-Cell. Also, this provided a pattern for displaying the 2-dimensional projective model of 4-variable actual problem.

**Keywords: 120-Cell, 600-Cell, 2-D projective model of 4-D rectangular coordinate system, 4-dimensional geometric object, 4-dimensional finite rotation group.**

## *Commemoration*

Harold Scott McDonald Coxeter (1907~2003), who died on the evening of March 31, 2003, aged 96, made fundamental contributions in the study of multi-dimensional geometric shapes and regarded as the greatest classical geometer of his generation.

Coxeter published in the geometrical field for 70 years, worked professionally at the

*University of Toronto for 60 years and wrote 12 books and more than 200 articles. He was best known for his work in hyper-dimensional geometries and **regular polytopes**--complicated geometric shapes of any number of dimensions that cannot be constructed in the real world but can be described mathematically and can sometime be drawn.*

*He wrote a book entitled **Regular Polytopes**[1] expounded specially shapes and constructions of the regular polytopes in $N$-dimensional space.*

*For commemorating Professor H.S.M.Coxeter, we plan to research deeply two most complicated regular polytopes **120-Cell** and **600-Cell** in $E^4$. In the meantime, some new results on shapes and data of these regular polytopes will be given.*

## 1. Introduction

Since 1850, Schläfli[2], Hilbert (Cohn-Vossen) and Coxeter investigated the basic situations of the regular polytopes in $E^4$ and higher, and made much progress. Their research might be helpful for people to understand the regular polytopes in $E^N$. Because of not having fundamentally established the 2-dimensional projective model of the *N*-dimensional rectangular coordinate systems, people cannot observe and investigate the true shapes of the regular polytopes in $E^N$.

For displaying the true shapes of two most complicated regular polytopes 120-Cell and 600-Cell, based on the mathematical theory, we plan to establish the 2-dimensional projective model of the 4-dimensional rectangular coordinate system, and adopt this projective model to deduce the transformation matrix, in order to calculate all vertex coordinates and determine the joint relationships between the adjacent vertices of two most complicated regular polytopes 120-Cell and 600-Cell, finally, to obtained some new data.

This research will be helpful for people to display and investigate the 4-variable actual problems in many fields.

## 2. To establish of 2-D projective model of 4-D rectangular coordinate system

Since *N*-dimensional space and the objects within it do exist, then, how to display them correctly and reasonably suggests an important topic to us.

For the 2-dimensional rectangular coordinate system, we can directly draw two axes which are in the orthogonal state on the plane to indicate it.

For the 3-dimensional rectangular coordinate system, we can use our right hand's thumb, index finger and middle finger which are in the each orthogonal state to stand for axes *X*, *Y* and *Z* respectively. We can also easily indicate the 2-dimensional projective model of 3-dimensional rectangular coordinate system with a plane figure, i.e. to draw axes *Y* and *Z* in orthogonal state and then to draw the positive semi-axis *X* to connect the positive semi-axes *Y* and *Z* at origin *O* and makes two angles of 135° with axes *Y* and *Z* respectively.

For the 4-dimensional rectangular coordinate system, to establish directly its model is a very difficult thing.

In order to overcome the obstacle of man's spatial thinking, we may as well make a further research for the 2-dimensional projective models of 4-dimensional objects.

Among all the regular polytopes in *N*-dimensional space, there is no more regular pattern which can directly embody the relations between all axes of the *N*-dimensional rectangular coordinate system than the hyper-cube.

Let's investigate the 4-dimensional cube. The following figure is the 2-dimensional projective model of 4-dimensional cube.

We may observe such fact: the 4-dimensional cube has 16 vertices, 32 edges, 24 squares and 8 cubes.

## 3. Using transformation matrix to simplify the 2-D plotting of 4-D geometric objects

In order to deduce the relevant transformation matrix, let's now establish the correspondence relationship between 4-dimensional geometric elements and their 2-dimensional projections. The following are the coordinates of 16 vertices of 4-dimensional cube.

1. (0,0,0,0)   5. (0,0,0,1)   9. (0,1,1,0)    13. (1,1,0,1)
2. (1,0,0,0)   6. (1,1,0,0)   10. (0,1,0,1)   14. (1,0,1,1)
3. (0,1,0,0)   7. (1,0,1,0)   11. (0,0,1,1)   15. (0,1,1,1)
4. (0,0,1,0)   8. (1,0,0,1)   12. (1,1,1,0)   16. (1,1,1,1)

So long as we demarcate these coordinates at the 2-dimensional projective model of 4-dimensional rectangular coordinate system (see Fig. 2. Notice: unit graduation of every coordinate axis equals to 1), then, we can get the 2-dimensional coordinates of above-mentioned vertices as follows:

1. $(0, 0)$;
2. $(-\frac{\sqrt{2}}{2}, -\frac{\sqrt{2}}{2})$
3. $(1, 0)$
4. $(0, 1)$
5. $(\frac{\sqrt{2}}{2}, -\frac{\sqrt{2}}{2})$;
6. $(1-\frac{\sqrt{2}}{2}, -\frac{\sqrt{2}}{2})$
7. $(-\frac{\sqrt{2}}{2}, 1-\frac{\sqrt{2}}{2})$
8. $(0, -\sqrt{2})$
9. $(1, 1)$;
10. $(1+\frac{\sqrt{2}}{2}, -\frac{\sqrt{2}}{2})$
11. $(\frac{\sqrt{2}}{2}, 1-\frac{\sqrt{2}}{2})$
12. $(1-\frac{\sqrt{2}}{2}, 1-\frac{\sqrt{2}}{2})$
13. $(1, -\sqrt{2})$
14. $(0, 1-\sqrt{2})$
15. $(1+\frac{\sqrt{2}}{2}, 1-\frac{\sqrt{2}}{2})$
16. $(1, 1-\sqrt{2})$

Because we have

$$(x^1, \ x^2, \ x^3, \ x^4) \begin{pmatrix} a_{11} & a_{12} \\ a_{21} & a_{22} \\ a_{31} & a_{32} \\ a_{41} & a_{42} \end{pmatrix} = (y^1, \ y^2),$$

therefore, so long as we substitut the coordinates of 16 vertices of 4-dimensional cube and their 2-dimensional coordinates into above formula, after calculating, we can obtain the transformation matrix of the geometric elements from 4-dimensional space project to 2-dimensional space:

$$A_2^4 = \begin{pmatrix} -\frac{\sqrt{2}}{2} & -\frac{\sqrt{2}}{2} \\ 1 & 0 \\ 0 & 1 \\ \frac{\sqrt{2}}{2} & -\frac{\sqrt{2}}{2} \end{pmatrix} \qquad (*)$$

where, the number of right-up of letter $A$ denotes original dimensional number of the geometric element, and the number of right-down denotes the projected dimensional number.

What should be noticed more is that from the 2-dimensional projective model of 4-dimensional cube, we discover that so long as we use following:

**Formula 1**. Let $A(x_1, x_2, ..., x_n)$ and $B(y_1, y_2, ..., y_n)$ be two points in $E^N$, then, the distance between $A$ and $B$ is

$$d = \sqrt{\sum_{i=1}^{n}(x_i - y_i)^2} \ ;$$

**Formula 2**. Let the coordinate of the vector $\boldsymbol{n}_1$ be $\{A_1, A_2, ..., A_n\}$ and the coordinate of the vector $\boldsymbol{n}_2$ be $\{B_1, B_2, ..., B_n\}$, then, the included angle between these two vectors is

$$\varphi = \cos^{-1}\left[\frac{\sum_{i=1}^{n} A_i B_i}{\sqrt{\sum_{i=1}^{n} A_i^2} \sqrt{\sum_{i=1}^{n} B_i^2}}\right].$$

We can know that the distance between any two most adjacent vertices of 4-dimensional cube are all equal to unit length 1. Undergoing such special projective transformation, we can find the property of unchanged distance between any two adjacent vertices, which is helpful to establish the 4-variable models by using the 2-dimensional projective model of 4-dimensional rectangular coordinate system, because it avoided the trouble of distortion coefficient of each axis.

## 4. To calculate vertex coordinates of the regular polytopes and determine the joint relationship between adjacent vertices

By the transformation matrix (*), first, we calculate the 2-dimensional projective coordinates of 4 adjacent vertices $(x_i^1, x_i^2, x_i^3, x_i^4)$ $(i = 1, 2, 3, 4)$ of the regular polytopes, and indicate in the 2-dimensional projective model of the 4-dimensional rectangular coordinate system, then, based on these vertices to list following equation group:

$$\begin{cases} (x^1 - x_1^1)^2 + (x^2 - x_1^2)^2 + (x^3 - x_1^3)^2 + (x^4 - x_1^4)^2 = a^2, \\ (x^1 - x_2^1)^2 + (x^2 - x_2^2)^2 + (x^3 - x_2^3)^2 + (x^4 - x_2^4)^2 = a^2, \\ (x^1 - x_3^1)^2 + (x^2 - x_3^2)^2 + (x^3 - x_3^3)^2 + (x^4 - x_3^4)^2 = a^2, \\ (x^1 - x_4^1)^2 + (x^2 - x_4^2)^2 + (x^3 - x_4^3)^2 + (x^4 - x_4^4)^2 = b^2. \end{cases}$$

where, $a$ is the distance between two most adjacent vertices (namely, $a$ is the length of an edge of the 4-dimensional regular polytopes), and $b$ is the distance between two rather adjacent vertices, certainly, $a < b$. After calculating, we can obtain the solution of this equation group as $(x_0^1,\ x_0^2,\ x_0^3,\ x_0^4)$.

Because we adopt the 2-dimensional projective model of 4-dimensional rectangular coordinate system, when we establish an equation group, in order to calculate the coordinate of next vertex, will feel have a certain direction. Otherwise, the people whom only have 3-dimensional spatial imagination, when they imagine the position and the shapes of 4-dimensional geometric objects, will be confused like a tangle of flax.

Also, according to the fixed distance $a$ between a pair of the adjacent vertices of the regular polytopes, we can obtain the joint relationship between all adjacent vertices of the regular polytopes.

**The table of 4-dimensional vertex coordinates (accuracy $10^{-9}$) of the 120-Cell and the 600-Cell and the table of the joint relationship between the adjacent vertices are listed below.**

## 5. To get new data of the 120-Cell and the 600-Cell

In addition, because we may draw the graphics of 4-dimensional geometric objects by using the 2-dimensional plotting system, therefore when we calculate concretely the date of 4-dimensional geometric objects (for examples, the vertex coordinate, the length of line segment, the angle between two line segments, the angle between two planes and the angle between two solids etc.) we will feel very convergence.

So long as we understand the volume formulas for solving the 3-dimensional geometric objects and the formulas for solving the angles between two intersect straight lines, two intersect planes and two intersect hyper-planes can be extended to 4-dimensional and higher, we will obtain much unknown data about the 120-Cell and the 600-Cell.

Now, we list the new data of the 120-Cell and the 600-Cell as follows:

| **Designation:** | **120-Cell** | | | **Schläfli Symbol: {5, 3, 3}** | | | |
|---|---|---|---|---|---|---|---|
| **I**. Composition | **II**. Feature | | | **III**. Numerical Value | | | |
| 600 vertices | 4 edges | meet | vertex | 1.The angle between two adjacent | edges | | 108° |
| 1200 edges | 3 pentagons | at a | edge | | pentagons | is | 116.5650507° |
| 720 pentagons | 2 dodecahedrons | (an) | pentagon | | dodecahedrons | | 144° |
| 120 dodecahedrons | | | | 2.Superficial Area | $919.5742838a^3$ | | |
| | | | | 3.Volume | $787.8569889a^4$ | | |

| **Designation:** | **600-Cell** | | | **Schläfli Symbol: {3, 3, 5}** | | | |
|---|---|---|---|---|---|---|---|
| **I.** Composition | **II.** Feature | | | **III**. Numerical Value | | | |
| 120 vertices | 12 edges | meet | vertex | 1. The angle between two adjacent | edges | | 60° |
| 720 edges | 5 triangles | at a | edge | | triangles | is | 70.52877936° |
| 1200 triangles | 2 tetrahedrons | (an) | triangle | | tetrahedrons | | 164.4775174° |
| 600 tetrahedrons | | | | 2.Superficial Area | $70.710678a^3$ | | |
| | | | | 3.Volume | $26.475425a^4$ | | |

# 6. To investigate the finite rotation groups of the 120-Cell and the 600-Cell

The pure mathematical theory has proved that the concrete geometric objects exist in *N*-dimensional space. We have discussed it in the reference [3, 4]. The following is an important theorem and some relevant formulas:

**Theorem ( Pole number theorem )** *Suppose that the order of finite rotation group of the regular polytopes in $E^N$ is $^N E(= o(^N G))$, and its 3-dimensional cell number is $\prod_{k=4}^{N} {^k v_{p_k}}$, as well as the cardinal of orbit is $^N v_{p_i}$, and the order of stability group of a point in orbit is $^N n_{p_i}$, then the following relationship holds:*

$$2(^N n - \prod_{k=4}^{N} {^k v_{p_k}}) = \sum_{i=1}^{N} {^N v_{p_i}} \cdot {^N n_{p_i}} \left[ 1 - \frac{^2 e_i}{p \cdot^{3-i} e_{3-i} \frac{2-i}{\{p,q\}}} \right].$$

The items $^N n$, $^N v_{p_i}$, $^N n_{p_i}$ of above relationship can be obtained from the following formulas:

**Formula 1(Order of finite rotation group)**
$$^N n = {^N v_{p_N}} \prod_{j=2}^{N-1} {^j e_j}, \qquad (i=1, 2,\cdots, N);$$

**Formula 2(Cardinal of orbit)**
$$^N v_{p_i} = \frac{^N v_{p_N} \cdot^{N-1} e_i}{^{N-i} e_{N-i} \frac{N-i-1}{\{\cdots,v,w\}}}, \qquad (i=1, 2,\cdots, N);$$

**Formula 3(Order of stability group)**
$$^N n_{p_i} = \frac{\prod_{j=2}^{N-1} {^j e_j} \cdot^{N-i} e_{N-i} \frac{N-i-1}{\{\cdots,v,w\}}}{^{N-1} e_i} \qquad (i = 1,2,\cdots, N);$$

**Formula 4(Unit element numbers of stability group)**
When $i=1,2,3$,
$$^N \alpha_{p_i} = \frac{^N n_{p_i} \cdot^2 e_i}{p \cdot^{3-i} e_{3-i} \frac{2-i}{\{p,q\}}};$$

when $i= 4, 5, \cdots, N$,

$$^N\alpha_{p_i} = {}^N n_{p_i}.$$

Now, we take the 120-Cell and the 600-Cell as the examples to verify conclusively the correctness of above theorem and formulas.

| Name | 120-Cell | | 600-Cell | |
|---|---|---|---|---|
| $^N n$ | 7200 | | 7200 | |
| $\prod_{k=4}^{N} {}^k v_{p_k}$ | 120 | | 600 | |
| $^N v_{p_i}$ | $i=1$ | 600 | $i=1$ | 120 |
| | $i=2$ | 1200 | $i=2$ | 720 |
| | $i=3$ | 720 | $i=3$ | 1200 |
| $^N n_{p_i}$ | $i=1$ | 12 | $i=1$ | 60 |
| | $i=2$ | 6 | $i=2$ | 10 |
| | $i=3$ | 10 | $i=3$ | 6 |
| $^N \alpha_{p_i}$ | $i=1$ | 4 | $i=1$ | 20 |
| | $i=2$ | 3 | $i=2$ | 5 |
| | $i=3$ | 2 | $i=3$ | 2 |
| **Contrast** | left | 14160 | left | 13200 |
| | right | 14160 | right | 13200 |
| **Result** | equals | | equals | |

This explains that the concrete objects exist in $N$-dimensional space, which enables the conception of $N$-dimensional space to break away from the abstract supposition so as to bring it back to reality again.

## Acknowledgements

I thank all those who helped in the preparation of this paper. In particular, I am grateful to Prof. Jiaoying Shi, Prof. Qunsheng Peng, Dr. Lizhuang Ma of Zhejiang University and Prof. Qizhao Ying , Associate Prof.

Tianhui Liu of Zhejiang Ocean University for their help and discussions. Also, I am grateful to Lecturer Jianmin Zhou for his help in computer modeling.

# The 4-Dimensional Coordinates of Vertices of the 120-Cell

| No | $x_1$ | $x_2$ | $x_3$ | $x_4$ |
|---|---|---|---|---|
| 1  | 0. | 0. | -0.8646837645 | -1.2228474870 |
| 2  | 0. | -0.2314036982 | 0.5659428717 | 1.3671850620 |
| 3  | 0. | -0.2314036982 | -0.8331448501 | 1.2228474870 |
| 4  | 0. | -0.2314036982 | 1.1003468280 | 0.9893043979 |
| 5  | 0. | 0.2314036982 | -1.4306266360 | 0.3778806546 |
| 6  | 0. | -0.3744190490 | -0.4833729193 | 1.3671850620 |
| 7  | 0. | 0.3744190490 | -1.4501187420 | 0. |
| 8  | 0. | -0.5174343998 | -0.9982847540 | 0.9893043979 |
| 9  | 0. | -0.5174343998 | 1.2654867320 | -0.6114237430 |
| 10 | 0. | -0.6058227472 | 0.0825699521 | -1.3671850620 |
| 11 | 0. | -0.6058227472 | 0.6169739091 | 1.2228474870 |
| 12 | 0. | -0.6058227472 | 0.9472537170 | 0.9893043979 |
| 13 | 0. | -0.6058227472 | -1.3165177690 | -.0.3778806546 |
| 14 | 0. | -0.7488380980 | -1.2970256470 | 0. |
| 15 | 0. | -0.7488380980 | 1.2970256470 | 0. |
| 16 | 0. | -0.7488380980 | -0.4323418823 | 1.2228474870 |
| 17 | 0. | -0.7488380980 | 0.4323418823 | -1.2228474870 |
| 18 | 0. | -0.8372264454 | -0.2161709413 | -1.2228474870 |
| 19 | 0. | -0.8372264454 | -0.7505748982 | 0.9893043979 |
| 20 | 0. | -0.8372264454 | -1.0808547060 | -0.6114237430 |
| 21 | 0. | -0.8372264454 | 1.1829167810 | -0.3778806546 |
| 22 | 0. | -1.0686301440 | 0.3497719302 | -0.9893043979 |
| 23 | 0. | -1.0686301440 | -1.0493157910 | 0. |
| 24 | 0. | -1.1232571470 | -0.0510310370 | -0.9893043979 |
| 25 | 0. | -1.1232571470 | -0.9157148022 | -0.3778806546 |
| 26 | 0. | -1.3546608450 | 0.5149118351 | 0.3778806546 |
| 27 | 0. | -1.3546608450 | 0.1846320260 | 0.6114237430 |
| 28 | 0. | -1.4430491930 | -0.1336009891 | 0.3778806546 |
| 29 | 0. | -1.4430491930 | 0.4008029672 | 0. |
| 30 | 0. | -1.4976761960 | 0. | 0. |
| 31 | -0.2022542486 | 0. | 0. | -1.4839566060 |
| 32 | 0.2022542486 | 0. | 0. | -1.4839566060 |
| 33 | -0.2022542486 | 0. | -1.3990877210 | -0.4946521989 |
| 34 | 0.2022542486 | 0. | -1.3990877210 | -0.4946521989 |
| 35 | -0.2022542486 | -0.3744190490 | 0.9157148022 | -1.1060759430 |
| 36 | 0.2022542486 | -0.3744190490 | 0.9157148022 | -1.1060759430 |
| 37 | -0.2022542486 | -0.3744190490 | -1.3480566840 | 0.4946521989 |
| 38 | -0.2022542486 | 0.3744190490 | 1.3480566840 | -0.4946521989 |
| 39 | -0.2022542486 | -0.6058227472 | -0.7821138131 | -1.1060759430 |
| 40 | 0.2022542486 | -0.6058227472 | -0.7821138131 | -1.1060759430 |
| 41 | -0.2022542486 | -0.9802417962 | 0.1336009891 | 1.1060759430 |
| 42 | 0.2022542486 | -0.9802417962 | 0.1336009891 | 1.1060759430 |
| 43 | -0.2022542486 | -0.9802417962 | 0.9982847540 | 0.4946521989 |
| 44 | 0.2022542486 | -0.9802417962 | 0.9982847540 | 0.4946521989 |
| 45 | -0.2022542486 | -1.2116454950 | -0.6995438612 | 0.4946521989 |
| 46 | -0.2022542486 | -1.2116454950 | 0.6995438612 | -0.4946521989 |
| 47 | 0.2022542486 | -1.2116454950 | 0.6995438612 | -0.4946521989 |
| 48 | 0.2022542486 | -1.2116454950 | -0.6995438612 | 0.4946521989 |
| 49 | -0.2022542486 | -1.3546608450 | -0.3497719302 | -0.4946521989 |
| 50 | 0.2022542486 | -1.3546608450 | -0.3497719302 | -0.4946521989 |
| 51 | -0.3272542486 | 0. | -0.8646837645 | 1.1782447260 |
| 52 | 0.3272542486 | 0. | -0.8646837645 | 1.1782447260 |
| 53 | -0.3272542486 | 0. | -1.3990877210 | 0.4224834156 |
| 54 | -0.3272542486 | 0. | 1.3990877210 | -0.4224834156 |
| 55 | -0.3272542486 | -0.2314036982 | -0.2987408932 | 1.4117878230 |
| 56 | -0.3272542486 | 0.2314036982 | 0.2987408932 | -1.4117878230 |
| 57 | -0.3272542486 | 0.1430153508 | -0.3497719302 | -1.4117878230 |
| 58 | 0.3272542486 | 0.1430153508 | -0.3497719302 | -1.4117878230 |
| 59 | -0.3272542486 | -0.1430153508 | 1.2144556950 | 0.8003640706 |
| 60 | -0.3272542486 | 0.1430153508 | -1.2144556950 | -0.8003640706 |
| 61 | -0.3272542486 | -0.2314036982 | -0.8331448501 | -1.1782447260 |

| 62 | 0.3272542486 | -0.2314036982 | -0.8331448501 | -1.1782447260 |
|---|---|---|---|---|
| 63 | -0.3272542486 | 0.2314036982 | -1.4306266350 | -0.1889403276 |
| 64 | 0.3272542486 | 0.2314036982 | -1.4306266350 | -0.1889403276 |
| 65 | -0.3272542486 | -0.3744190490 | 0.0510310370 | -1.4117878230 |
| 66 | 0.3272542486 | -0.3744190490 | 0.0510310370 | -1.4117878230 |
| 67 | -0.3272542486 | -0.3744190490 | -1.3480566840 | -0.4224834156 |
| 68 | 0.3272542486 | -0.3744190490 | -1.3480566840 | -0.4224834156 |
| 69 | -0.3272542486 | 0.4628073964 | 1.1318857430 | -0.8003640706 |
| 70 | 0.3272542486 | -0.4628073964 | 1.1318857430 | -0.8003540706 |
| 71 | -0.3272542486 | -0.4628073964 | -1.1318857430 | 0.8003640706 |
| 72 | -0.3272542486 | -0.4628073964 | 1.1318857430 | -0.8003640706 |
| 73 | -0.3272542486 | -0.6058227472 | 0.6169739091 | -1.1782447260 |
| 74 | 0.3272542486 | -0.6058227472 | 0.6169739091 | -1.1782447260 |
| 75 | -0.3272542486 | -0.6058227472 | -1.3165177690 | 0.1889403276 |
| 76 | 0.3272542486 | -0.6058227472 | -1.3165177690 | 0.1889403276 |
| 77 | -0.3272542486 | -0.7488380980 | -0.4323418823 | -1.1782447260 |
| 78 | -0.3272542486 | -0.7488380980 | 0.4323418823 | 1.1782447260 |
| 79 | 0.3272542486 | -0.7488380980 | -0.4323418823 | -1.1782447260 |
| 80 | 0.3272542486 | -0.7488380980 | 0.4323418823 | 1.1782447260 |
| 81 | -0.3272541486 | -0.7488380980 | -0.9667458392 | -0.8003640706 |
| 82 | -0.3272542486 | -0.7488380980 | 0.9667458392 | 0.8003640706 |
| 83 | 0.3272542486 | -0.7488380980 | -0.9667458392 | -0.8003640706 |
| 84 | 0.3272542486 | -0.7488380980 | 0.9667458392 | 0.8003640706 |
| 85 | -0.3272542486 | -0.8372264454 | 1.1829167810 | 0.1889403276 |
| 86 | 0.3272542486 | -0.8372264454 | 1.1829167810 | 0.1889403276 |
| 87 | 0.3272542486 | -0.8373264454 | -0.2161709413 | 1.1782447260 |
| 88 | -0.3272542486 | -0.8372264454 | -0.2161709413 | 1.1782447260 |
| 89 | -0.3272542486 | -0.9802417962 | 0.9982847540 | -0.4224834156 |
| 90 | 0.3272542486 | -0.9802417962 | 0.9982847540 | -0.4224834156 |
| 91 | -0.3272542486 | -0.9802417962 | -0.7310827761 | 0.8003640706 |
| 92 | 0.3272542486 | -0.9802417962 | -0.7310827761 | 0.8003640706 |
| 93 | -0.3272542486 | -1.1232571470 | 0.4833729193 | -0.8003640706 |
| 94 | 0.3272542486 | -1.1232571470 | 0.4833729193 | -0.8003640706 |
| 95 | -0.3272542486 | -1.1232571470 | -0.9157148022 | 0.1889403276 |
| 96 | 0.3272542486 | -1.1232571470 | -0.9157148022 | 0.1889403276 |
| 97 | -0.3272542486 | -1.2116454950 | -0.6995438612 | -0.4224834156 |
| 98 | -0.3272542486 | -1.2116454950 | 0.6995438612 | 0.4224834156 |
| 99 | 0.3272542486 | -1.2116454950 | -0.6995438612 | -0.4224834156 |
| 100 | 0.3272542486 | -1.2116454950 | 0.6995438612 | 0.4224834156 |
| 101 | -0.3272542486 | -1.2116454950 | -0.1651399043 | -0.8003640706 |
| 102 | -0.3272542486 | -1.2116454950 | 0.1651399043 | 0.8003640706 |
| 103 | 0.3272542486 | -1.2116454950 | -0.1651399043 | -0.8003640706 |
| 104 | 0.3272542486 | -1.2116454950 | 0.1651399043 | 0.8003640706 |
| 105 | -0.3272542486 | -1.3546608450 | -0.3497719302 | 0.4224834156 |
| 106 | 0.3272542486 | -1.3546608450 | -0.3497719302 | 0.4224834156 |
| 107 | -0.3272542486 | -1.3546608450 | 0.5149118351 | -0.1889403276 |
| 108 | 0.3272542486 | -1.3546608450 | 0.5149118351 | -0.1889403276 |
| 109 | -0.3272542486 | -1.4430491930 | -0.1336009891 | -0.1889403276 |
| 110 | 0.3272542486 | -1.4430491930 | -0.1336009891 | -0.1889403276 |
| 111 | -0.5295084972 | 0. | -0.5344039562 | 1.2950162700 |
| 112 | -0.5295084972 | 0. | 0.5344039562 | -1.2950162700 |
| 113 | -0.5295084972 | 0. | -1.3990877210 | 0.0721687833 |
| 114 | 0.5295084972 | 0. | -1.3990877210 | 0.0721687833 |
| 115 | -0.5295084972 | -0.1430153508 | -0.5149118351 | -1.2950162700 |
| 116 | 0.5295084972 | -0.1430153508 | -0.5149118351 | -1.2950162700 |
| 117 | -0.5295084972 | -0.1430153508 | -1.0493157910 | -0.9171356145 |
| 118 | 0.5295084972 | -0.1430153508 | -1.0493157910 | -0.9171356145 |
| 119 | -0.5295084972 | -0.1430153508 | 1.2144556950 | -0.6835925263 |
| 120 | -0.5295084972 | 0.1430153508 | -1.2144556950 | 0.6835925263 |
| 121 | -0.5295084972 | -0.1430153508 | -1.0493157910 | 0.9171356145 |
| 122 | -0.5295084972 | 0.1430153508 | 1.0493157910 | -0.9171356145 |
| 123 | -0.5295084972 | -0.3744190490 | 0.3813108453 | -1.2950162700 |
| 124 | -0.5295084972 | 0.3744190490 | -0.3813108453 | 1.2950162700 |
| 125 | -0.5295084972 | -0.3744190490 | -1.3480566840 | -0.0721687833 |

| 126 | 0.5295084972 | -0.3744190490 | -1.3480566840 | -0.0721687833 |
|---|---|---|---|---|
| 127 | -0.5295084972 | -0.4628073964 | -0.2672019781 | -1.2950162700 |
| 128 | 0.5295084972 | -0.4628073964 | -0.2672019781 | -1.2950162700 |
| 129 | 0.5295084972 | -0.4628073964 | 0.2672019781 | 1.2950162700 |
| 130 | 0.5295084972 | 0.4628073964 | -0.2672019781 | -1.2950162700 |
| 131 | -0.5295084972 | -0.4628073964 | -1.1318857430 | -0.6835925263 |
| 132 | -0.5295084972 | -0.4628073964 | 1.1318857430 | 0.6835925263 |
| 133 | -0.5295084972 | 0.4628073964 | -1.1318857430 | -0.6835925263 |
| 134 | 0.5295084972 | -0.4628073964 | -1.1318857430 | -0.6835925263 |
| 135 | -0.5295084972 | -0.5174343998 | -0.1336009891 | 1.2950162700 |
| 136 | -0.5295084972 | 0.5174343998 | 0.1336009891 | -1.2950162700 |
| 137 | -0.5295084972 | -0.5174343998 | 1.2654867320 | 0.3057118717 |
| 138 | 0.5295084972 | -0.5174343998 | 1.2654867320 | 0.3057118717 |
| 139 | -0.5295084972 | -0.7488380980 | -0.9667458392 | 0.6835925263 |
| 140 | -0.5295084972 | -0.7488380980 | 0.9667458392 | -0.6835925263 |
| 141 | 0.5295084972 | -0.7488380980 | -0.9667458392 | 0.6835925263 |
| 142 | 0.5295084972 | -0.7488380980 | 0.9667458392 | -0.6835925263 |
| 143 | -0.5295084972 | -0.8372264454 | 0.6485128238 | -0.9171356154 |
| 144 | 0.5295084972 | -0.8372264454 | 0.6485128238 | -0.9171356154 |
| 145 | 0.5295084972 | -0.8372264454 | 0.6485128238 | 0.9171356154 |
| 146 | -0.5295084972 | -0.8372264454 | 0.6485128238 | 0.9171356154 |
| 147 | -0.5295084972 | -0.8372264454 | -1.0808547060 | 0.3057118717 |
| 148 | 0.5295084972 | -0.8372264454 | -1.0808547060 | 0.3057118717 |
| 149 | -0.5295084972 | -0.9802417962 | -0.4008029672 | 0.9171356154 |
| 150 | -0.5295084972 | -0.9802417962 | -0.4008029672 | -0.9171356154 |
| 151 | -0.5295084972 | -0.9802417962 | -0.7310827761 | -0.6835925263 |
| 152 | -0.5295084972 | 0.9802417962 | 0.7310827761 | 0.6835925263 |
| 153 | -0.5295084972 | 0.9802417962 | 0.4008029672 | 0.9171356154 |
| 154 | 0.5295084972 | -0.9802417962 | -0.4008029672 | 0.9171356154 |
| 155 | -0.5295084972 | -0.9802417962 | 0.9982847540 | -0.0721687833 |
| 156 | 0.5295084972 | -0.9802417962 | 0.9982847540 | -0.0721687833 |
| 157 | -0.5295084972 | -1.1232571470 | 0.4833729193 | 0.6835925263 |
| 158 | 0.5295084972 | -1.1232571470 | 0.4833729193 | 0.6835925263 |
| 159 | -0.5295084972 | -1.2116454950 | 0.1651399043 | -0.6835925263 |
| 160 | 0.5295084972 | -1.2116454950 | 0.1651399043 | -0.6835925263 |
| 161 | -0.5295084972 | -1.2116454950 | -0.6995438612 | -0.0721687833 |
| 162 | -0.5295084972 | -1.2116454950 | 0.6995438612 | 0.0721687833 |
| 163 | -0.5295084972 | -1.2116454950 | -0.1651399043 | 0.6835925263 |
| 164 | 0.5295084972 | -1.2116454950 | -0.1651399043 | 0.6835925263 |
| 165 | 0.5295084972 | -1.2116454950 | 0.6995438612 | 0.0721687833 |
| 166 | 0.5295084972 | -1.2116454950 | -0.6995438612 | -0.0721687833 |
| 167 | -0.5295084972 | -1.3546608450 | 0.1846320260 | -0.3057118717 |
| 168 | 0.5295084972 | -1.3546608450 | 0.1846320260 | -0.3057118717 |
| 169 | -0.5295084972 | -1.3546608450 | -0.3497719302 | 0.0721687833 |
| 170 | 0.5295084972 | -1.3546608450 | -0.3497719302 | 0.0721687833 |
| 171 | -0.8567627458 | 0. | -0.5344039562 | -1.1060759430 |
| 172 | -0.8567627458 | 0. | 0.5344039562 | 1.1060759430 |
| 173 | -0.8567627458 | 0. | -0.8646837645 | -0.8725328536 |
| 174 | -0.8567627458 | 0. | 0.8646837645 | 0.8725328536 |
| 175 | -0.8567627458 | -0.1430153508 | -0.5149118351 | 1.1060759430 |
| 176 | -0.8567627458 | 0.1430153508 | 0.5149118351 | -1.1060759430 |
| 177 | -0.8567627458 | -0.1430153508 | 1.2144556950 | -0.1167715443 |
| 178 | -0.8567627458 | 0.1430153508 | -1.2144556950 | 0.1167715443 |
| 179 | -0.8567627458 | -0.2314036982 | 1.1003468280 | -0.4946521989 |
| 180 | -0.8567627458 | 0.2314036982 | -1.1003468280 | 0.4946521989 |
| 181 | -0.8567627458 | -0.2314036982 | -0.8331448501 | 0.8725328536 |
| 182 | 0.8567627458 | -0.2314036982 | -0.8331448501 | 0.8725328536 |
| 183 | -0.8567627458 | -0.3744190490 | 0.9157148022 | 0.7281952873 |
| 184 | -0.8567627458 | 0.3744190490 | -0.9157148022 | -0.7281952873 |
| 185 | -0.8567627458 | -0.3744190490 | 0.3813108453 | 1.1060759430 |
| 186 | 0.8567627458 | -0.3744190490 | 0.3813108453 | 1.1060759430 |
| 187 | -0.8567627458 | -0.4628073964 | 0.2672019781 | -1.1060759430 |
| 188 | -0.8567627458 | -0.4628073964 | -0.2672019781 | 1.1060759430 |
| 189 | -0.8567627458 | 0.4628073964 | 0.2672019781 | -1.1060759430 |
| 190 | 0.8567627458 | -0.4628073964 | 0.2672019781 | -1.1060759430 |
| 191 | -0.8567627458 | -0.4628073964 | -1.1318857430 | -0.1167715443 |
| 192 | -0.8567627458 | 0.4628073964 | -1.1318857430 | -0.1167715443 |
| 193 | -0.8567627458 | -0.4628073964 | 1.1318857430 | 0.1167715443 |
| 194 | 0.8567627458 | -0.4628073964 | -1.1318857430 | -0.1167715443 |
| 195 | -0.8567627458 | -0.5174343998 | -0.1336009891 | -1.1060759430 |
| 196 | -0.8567627458 | 0.5174343998 | 0.1336009891 | 1.1060759430 |
| 197 | -0.8567627458 | -0.5174343998 | -0.9982847540 | -0.4946521989 |
| 198 | -0.8567627458 | 0.5174343998 | 0.9982847540 | 0.4946521989 |
| 199 | -0.8567627458 | -0.6058227472 | -0.7821138131 | 0.7281952873 |
| 200 | 0.8567627458 | -0.6058227472 | -0.7821138131 | 0.7281952873 |
| 201 | -0.8567627458 | -0.6058227472 | 0.6169739091 | 0.8725328536 |
| 202 | 0.8567627458 | -0.6058227472 | 0.6169739091 | 0.8725328536 |
| 203 | -0.8567627458 | 0.6058227472 | -0.9472537170 | 0.4946521989 |
| 204 | -0.8567627458 | -0.6058227472 | 0.9472537170 | -0.4946521989 |
| 205 | -0.8567627458 | -0.7488380980 | -0.4323418823 | 0.8725328536 |
| 206 | -0.8567627458 | -0.7488380980 | 0.4323418823 | -0.8725328536 |
| 207 | 0.8567627458 | -0.7488380980 | -0.4323418823 | 0.8725328536 |
| 208 | 0.8567627458 | -0.7488380980 | 0.4323418823 | -0.8725328536 |
| 209 | -0.8567627458 | -0.7488380980 | -0.9667458392 | 0.1167715443 |
| 210 | -0.8567627458 | -0.7488380980 | 0.9667458392 | -0.11677154430 |
| 211 | 0.8567627458 | -0.7488380980 | -0.9667458392 | 0.1167715443 |
| 212 | 0.8567627458 | -0.7488380980 | 0.9667458392 | -0.1167715443 |
| 213 | -0.8567627458 | -0.8372264454 | -0.2161709413 | -0.8725328536 |
| 214 | 0.8567627458 | -0.8372264454 | -0.2161709413 | -0.8725328536 |
| 215 | -0.8567627458 | -0.8372264454 | -0.7505748982 | -0.4946521989 |
| 216 | 0.8567627458 | -0.8372264454 | -0.7505748982 | -0.4946521989 |
| 217 | -0.8567627458 | -0.9802417962 | 0.1336009891 | -0.7281952873 |
| 218 | 0.8567627458 | -0.9802417962 | 0.1336009891 | -0.7281952873 |
| 219 | -0.8567627458 | -0.9802417962 | -0.7310827761 | -0.1167715443 |
| 220 | 0.8567627458 | -0.9802417962 | -0.7310827761 | -0.1167715443 |
| 221 | -0.8567627458 | -1.0686301440 | 0.3497719302 | 0.4946521989 |

| 222 | 0.8567627458 | -1.0686301440 | 0.3497719302 | 0.4946521989 |
|---|---|---|---|---|
| 223 | -0.8567627458 | -1.1232571470 | 0.4833729193 | 0.1167715443 |
| 224 | 0.8567627458 | -1.1232571470 | 0.4833729193 | 0.1167715443 |
| 225 | -0.8567627458 | -1.1232571470 | -0.0510310370 | 0.4946521989 |
| 226 | 0.8567627458 | -1.1232571470 | -0.0510310370 | 0.4946521989 |
| 227 | -0.8567627485 | -1.2116454950 | -0.1651399043 | 0.1167715443 |
| 228 | -0.8567627458 | -1.2116454950 | 0.1651399043 | -0.1167715443 |
| 229 | -0.8567627458 | 1.2116454950 | 0.1651399043 | -0.1167715443 |
| 230 | 0.8567627458 | -1.2116454950 | 0.1651399043 | -0.1167715443 |
| 231 | -1.0590169944 | 0. | -0.8646837645 | 0.6114237430 |
| 232 | -1.0590169944 | 0. | 0.8646837645 | -0.6114237430 |
| 233 | -1.0590169944 | -0.1430153508 | 0.3497719302 | -0.9893043979 |
| 234 | -1.0590169944 | 0.1430153508 | -0.3497719302 | 0.9893043979 |
| 235 | -1.0590169944 | -0.1430153508 | -1.0493157910 | 0. |
| 236 | -1.0590169944 | 0.1430153508 | 1.0493157910 | 0. |
| 237 | -1.0590169944 | -0.2314036982 | -0.2987408932 | -0.9893043979 |
| 238 | -1.0590169944 | 0.2314036982 | 0.2987408932 | 0.9893043979 |
| 239 | -1.0590169944 | -0.2314036982 | -0.8331448501 | -0.6114237430 |
| 240 | -1.0590169944 | 0.2314036982 | 0.8331448501 | 0.6114237430 |
| 241 | -1.0590169944 | -0.3744190490 | 0.0510310370 | 0.9893043979 |
| 242 | -1.0590169944 | 0.3744190490 | -0.0510310370 | -0.9893043979 |
| 243 | -1.0590169944 | -0.3744190490 | 0.9157148022 | 0.3778806546 |
| 244 | -1.0590169944 | 0.3744190490 | -0.9157148022 | -0.3778806546 |
| 245 | -1.0590169944 | -0.6058227472 | 0.6169739091 | -0.6114237430 |
| 246 | 1.0590169944 | -0.6058227472 | 0.6169739091 | -0.6114237430 |
| 247 | -1.0590169944 | -0.6058227472 | -0.7821138131 | 0.3778806546 |
| 248 | 1.0590169944 | -0.6058227472 | -0.7821138131 | 0.3778806546 |
| 249 | -1.0590169944 | -0.7488380980 | -0.4323418823 | -0.6114237430 |
| 250 | -1.0590169944 | -0.7488380980 | 0.4323418823 | 0.6114237430 |
| 251 | -1.0590169944 | 0.7488380980 | 0.4323418823 | 0.6114237430 |
| 252 | -1.0590169944 | 0.7488380980 | -0.4323418823 | -0.6114237430 |
| 253 | -1.0590169944 | -0.8372264454 | -0.2161709413 | 0.6114237430 |

| 254 | 1.0590169944 | -0.8372264454 | -0.2161709413 | 0.6114237430 |
| 255 | -1.0590169944 | -0.8372264454 | 0.6485128238 | 0. |
| 256 | 1.0590169944 | -0.8372264454 | 0.6485128238 | 0. |
| 257 | -1.0590169944 | -0.9802417962 | 0.1336009891 | -0.3778806546 |
| 258 | 1.0590169944 | -0.9802417962 | 0.1336009891 | -0.3778806546 |
| 259 | -1.0590169944 | -0.9802417962 | -0.4008029672 | 0. |
| 260 | 1.0590169944 | -0.9802417962 | -0.4008029672 | 0. |
| 261 | -1.1840169944 | 0. | 0. | -0.9171356154 |
| 262 | -1.1840169944 | 0. | 0. | 0.9171356154 |
| 263 | -1.1840169944 | 0. | -0.8646837645 | -0.3057118717 |
| 264 | -1.1840169944 | 0. | 0.8646837645 | 0.3057118717 |
| 265 | -1.1840169944 | -0.2314036982 | -0.8331448501 | 0.3057118717 |
| 266 | -1.1840169944 | 0.2314036982 | 0.8331448501 | -0.3057118717 |
| 267 | -1.1840169944 | -0.2314036982 | 0.5659428717 | -0.6835925263 |
| 268 | -1.1840169944 | 0.2314036982 | -0.5659428717 | 0.6835925263 |
| 269 | -1.1840169944 | -0.3744190490 | -0.4833729193 | -0.6835925263 |
| 270 | -1.1840169944 | 0.3744190490 | 0.4833729193 | 0.6835925263 |
| 271 | -1.1840169944 | -0.6058227472 | 0.0825688521 | 0.6835925263 |
| 272 | -1.1840169944 | 0.6058227472 | -0.0825688521 | -0.6835925263 |
| 273 | -1.1840169944 | -0.6058227472 | 0.6169739091 | 0.3057118717 |
| 274 | 1.1840169944 | -0.6058227472 | 0.6169739091 | 0.3057118717 |
| 275 | -1.1840169944 | -0.7488380980 | 0.4323418823 | -0.3057118717 |
| 276 | -1.1840169944 | -0.7488380980 | -0.4323418823 | 0.3057118717 |
| 277 | -1.1840169944 | 0.7488380980 | 0.4323418823 | -0.3057118717 |
| 278 | 1.1840169944 | -0.7488380980 | 0.4323418823 | -0.3057118717 |
| 279 | -1.1840169944 | -0.8372264454 | -0.2161709413 | -0.3057118717 |
| 280 | -1.1840169944 | 0.8372264454 | 0.2161709413 | 0.3057118717 |
| 281 | -1.3862712430 | 0. | 0. | -0.5668209828 |
| 282 | -1.3862712430 | 0. | 0. | 0.5668209828 |
| 283 | -1.3862712430 | 0. | -0.5344039562 | -0.1889403276 |
| 284 | -1.3862712430 | 0. | 0.5344039562 | 0.1889403276 |
| 285 | -1.3862712430 | -0.1430153508 | 0.3497719302 | -0.4224834156 |

| 286 | -1.3862712430 | 0.1430153508 | -0.3497719302 | 0.4224834156 |
|---|---|---|---|---|
| 287 | -1.3862712430 | -0.1430153508 | -0.5149118351 | 0.1889403276 |
| 288 | -1.3862712430 | 0.1430153508 | 0.5149118351 | -0.1889403276 |
| 289 | -1.3862712430 | -0.2314036982 | -0.2987408932 | -0.4224834156 |
| 290 | -1.3862712430 | 0.2314036982 | 0.2987408932 | 0.4224834156 |
| 291 | -1.3862712430 | -0.3744190490 | 0.0510310370 | 0.4224834156 |
| 292 | -1.3862712430 | 0.3744190490 | -0.0510310370 | -0.4224834156 |
| 293 | -1.3862712430 | -0.3744190490 | 0.3813108453 | 0.1889403276 |
| 294 | -1.3862712430 | 0.3744190490 | -0.3813108453 | -0.1889403276 |
| 295 | -1.3862712430 | -0.4628073964 | -0.2672019781 | 0.1889403276 |
| 296 | -1.3862712430 | -0.4628073964 | 0.2672019781 | -0.1889403276 |
| 297 | -1.3862712430 | 0.4628073964 | -0.2672019781 | 0.1889403276 |
| 298 | -1.3862712430 | 0.4628073964 | 0.2672019871 | -0.1889403276 |
| 299 | -1.3862712430 | -0.5174343998 | -0.1336009891 | -0.1889403276 |
| 300 | -1.3862712430 | 0.5174343998 | 0.1336009891 | 0.1889403276 |
| 301 | 0. | 0. | 0.8646837645 | 1.2228474870 |
| 302 | 0. | 0.2314036982 | -0.5659428717 | -1.3671850620 |
| 303 | 0. | 0.2314036982 | 0.8331448501 | -1.2228474870 |
| 304 | 0. | 0.2314036982 | -1.1003468280 | -0.9893043979 |
| 305 | 0. | -0.2314036982 | 1.4306266360 | -0.3778806546 |
| 306 | 0. | 0.3744190490 | 0.4833729193 | -1.3671850620 |
| 307 | 0. | -0.3744190490 | 1.4501187420 | 0. |
| 308 | 0. | 0.5174343998 | 0.9982847540 | -0.9893043979 |
| 309 | 0. | 0.5174343998 | -1.2654867320 | 0.6114237430 |
| 310 | 0. | 0.6058227472 | -0.0825699521 | 1.3671850620 |
| 311 | 0. | 0.6058227472 | -0.6169739091 | -1.2228474870 |
| 312 | 0. | 0.6058227472 | -0.9472537170 | -0.9893043979 |
| 313 | 0. | 0.6058227472 | 1.3165177690 | 0.3778806546 |
| 314 | 0. | 0.7488380980 | 1.2970256470 | 0. |
| 315 | 0. | 0.7488380980 | -1.2970256470 | 0. |
| 316 | 0. | 0.7488380980 | 0.4323418823 | -1.2228474870 |
| 317 | 0. | 0.7488380980 | -0.4323418823 | 1.2228474870 |

| 318 | 0. | 0.8372264454 | 0.2161709413 | 1.2228474870 |
|---|---|---|---|---|
| 319 | 0. | 0.8372264454 | 0.7505748982 | -0.9893043979 |
| 320 | 0. | 0.8372264454 | 1.0808547060 | 0.6114237430 |
| 321 | 0. | 0.8372264454 | -1.1829167810 | 0.3778806546 |
| 322 | 0. | 1.0686301440 | -0.3497719302 | 0.9893043979 |
| 323 | 0. | 1.0686301440 | 1.0493157910 | 0. |
| 324 | 0. | 1.1232571470 | 0.0510310370 | 0.9893043979 |
| 325 | 0. | 1.1232571470 | 0.9157148022 | 0.3778806546 |
| 326 | 0. | 1.3546608450 | -0.5149118351 | -0.3778806546 |
| 327 | 0. | 1.3546608450 | -0.5149118351 | -0.3778806546 |
| 328 | 0. | 1.4430491930 | 0.1336009891 | -0.3778806546 |
| 329 | 0. | 1.4430491930 | -0.4008029672 | 0. |
| 330 | 0. | 1.4976761960 | 0. | 0. |
| 331 | 0.2022542486 | 0. | 0. | 1.4839566060 |
| 332 | -0.2022542486 | 0. | 0. | 1.4839566060 |
| 333 | 0.2022542486 | 0. | 1.3990877210 | 0.4946521989 |
| 334 | -0.2022542486 | 0. | 1.3990877210 | 0.4946521989 |
| 335 | 0.2022542486 | 0.3744190490 | -0.9157148022 | 1.1060759430 |
| 336 | -0.2022542486 | 0.3744190490 | -0.9157148022 | 1.1060759430 |
| 337 | 0.2022542486 | 0.3744190490 | 1.3480566840 | -0.4946521989 |
| 338 | 0.2022542486 | -0.3744190490 | -1.3480566840 | 0.4946521989 |
| 339 | 0.2022542486 | 0.6058227472 | 0.7821138131 | 1.1060759430 |
| 340 | -0.2022542486 | 0.6058227472 | 0.7821138131 | 1.1060759430 |
| 341 | 0.2022542486 | 0.9802417962 | -0.1336009891 | -1.1060759430 |
| 342 | - 0.2022542486 | 0.9802417962 | -0.1336009891 | -1.1060759430 |
| 343 | 0.2022542486 | 0.9802417962 | -0.9982847540 | -0.4946521989 |
| 344 | -0.2022542486 | 0.9802417962 | -0.9982847540 | -0.4946521989 |
| 345 | 0.2022542486 | 1.2116454950 | 0.6995438612 | -0.4946521989 |
| 346 | 0.2022542486 | 1.2116454950 | -0.6995438612 | 0.4946521989 |
| 347 | -0.2022542486 | 1.2116454950 | -0.6995438612 | 0.4946521989 |
| 348 | -0.2022542486 | 1.2116454950 | 0.6995438612 | -0.4946521989 |
| 349 | 0.2022542486 | 1.3546608450 | 0.3497719302 | 0.4946521989 |

| 350 | -0.2022542486 | 1.3546608450 | 0.3497719302 | 0.4946521989 |
|---|---|---|---|---|
| 351 | 0.3272542486 | 0. | 0.8646837645 | -1.1782447260 |
| 352 | -0.3272542486 | 0. | 0.8646837645 | -1.1782447260 |
| 353 | 0.3272542486 | 0. | 1.3990877210 | -0.4224834156 |
| 354 | 0.3272542486 | 0. | -1.3990877210 | 0.4224834156 |
| 355 | 0.3272542486 | 0.2314036982 | 0.2987408932 | -1.4117878230 |
| 356 | 0.3272542486 | -0.2314036982 | -0.2987408932 | 1.4117878230 |
| 357 | 0.3272542486 | -0.1430153508 | 0.3497719302 | 1.4117878230 |
| 358 | -0.3272542486 | -0.1430153508 | 0.3497719302 | 1.4117878230 |
| 359 | 0.3272542486 | 0.1430153508 | -1.2144556950 | -0.8003640706 |
| 360 | 0.3272542486 | -0.1430153508 | 1.2144556950 | 0.8003640706 |
| 361 | 0.3272542486 | 0.2314036982 | 0.8331448501 | 1.1782447260 |
| 362 | -0.3272542486 | 0.2314036982 | 0.8331448501 | 1.1782447260 |
| 363 | 0.3272542486 | -0.2314036982 | 1.4306266350 | 0.1889403276 |
| 364 | -0.3272542486 | -0.2314036982 | 1.4306266350 | 0.1889403276 |
| 365 | 0.3272542486 | 0.3744190490 | -0.0510310370 | 1.4117878230 |
| 366 | -0.3272542486 | 0.3744190490 | -0.0510310370 | 1.4117878230 |
| 367 | 0.3272542486 | 0.3744190490 | 1.3480566840 | 0.4224834156 |
| 368 | -0.3272542486 | 0.3744190490 | 1.3480566840 | 0.4224834156 |
| 369 | 0.3272542486 | -0.4628073964 | -1.1318857430 | 0.8003640706 |
| 370 | -0.3272542486 | 0.4628073964 | -1.1318857430 | 0.8003540706 |
| 371 | 0.3272542486 | 0.4628073964 | 1.1318857430 | -0.8003640706 |
| 372 | 0.3272542486 | 0.4628073964 | -1.1318857430 | 0.8003640706 |
| 373 | 0.3272542486 | 0.6058227472 | -0.6169739091 | 1.1782447260 |
| 374 | -0.3272542486 | 0.6058227472 | -0.6169739091 | 1.1782447260 |
| 375 | 0.3272542486 | 0.6058227472 | 1.3165177690 | -0.1889403276 |
| 376 | -0.3272542486 | 0.6058227472 | 1.3165177690 | -0.1889403276 |
| 377 | 0.3272542486 | 0.7488380980 | 0.4323418823 | 1.1782447260 |
| 378 | 0.3272542486 | 0.7488380980 | -0.4323418823 | -1.1782447260 |
| 379 | -0.3272542486 | 0.7488380980 | 0.4323418823 | 1.1782447260 |
| 380 | -0.3272542486 | 0.7488380980 | -0.4323418823 | -1.1782447260 |
| 381 | 0.3272541486 | 0.7488380980 | 0.9667458392 | 0.8003640706 |
| 382 | 0.3272542486 | 0.7488380980 | -0.9667458392 | -0.8003640706 |
| 383 | -0.3272542486 | 0.7488380980 | 0.9667458392 | 0.8003640706 |
| 384 | -0.3272542486 | 0.7488380980 | -0.9667458392 | -0.8003640706 |
| 385 | 0.3272542486 | 0.8372264454 | -1.1829167810 | -0.1889403276 |
| 386 | -0.3272542486 | 0.8372264454 | -1.1829167810 | -0.1889403276 |
| 387 | -0.3272542486 | 0.8373264454 | 0.2161709413 | -1.1782447260 |
| 388 | 0.3272542486 | 0.8372264454 | 0.2161709413 | -1.1782447260 |
| 389 | 0.3272542486 | 0.9802417962 | -0.9982847540 | 0.4224834156 |
| 390 | -0.3272542486 | 0.9802417962 | -0.9982847540 | 0.4224834156 |
| 391 | 0.3272542486 | 0.9802417962 | 0.7310827761 | -0.8003640706 |
| 392 | -0.3272542486 | 0.9802417962 | 0.7310827761 | -0.8003640706 |
| 393 | 0.3272542486 | 1.1232571470 | -0.4833729193 | 0.8003640706 |
| 394 | -0.3272542486 | 1.1232571470 | -0.4833729193 | 0.8003640706 |
| 395 | 0.3272542486 | 1.1232571470 | 0.9157148022 | -0.1889403276 |
| 396 | -0.3272542486 | 1.1232571470 | 0.9157148022 | -0.1889403276 |
| 397 | 0.3272542486 | 1.2116454950 | 0.6995438612 | 0.4224834156 |
| 398 | 0.3272542486 | 1.2116454950 | -0.6995438612 | -0.4224834156 |
| 399 | -0.3272542486 | 1.2116454950 | 0.6995438612 | 0.4224834156 |
| 400 | -0.3272542486 | 1.2116454950 | -0.6995438612 | -0.4224834156 |
| 401 | 0.3272542486 | 1.2116454950 | 0.1651399043 | 0.8003640706 |
| 402 | 0.3272542486 | 1.2116454950 | -0.1651399043 | -0.8003640706 |
| 403 | -0.3272542486 | 1.2116454950 | 0.1651399043 | 0.8003640706 |
| 404 | -0.3272542486 | 1.2116454950 | -0.1651399043 | -0.8003640706 |
| 405 | 0.3272542486 | 1.3546608450 | 0.3497719302 | -0.4224834156 |
| 406 | -0.3272542486 | 1.3546608450 | 0.3497719302 | -0.4224834156 |
| 407 | 0.3272542486 | 1.3546608450 | -0.5149118351 | 0.1889403276 |
| 408 | -0.3272542486 | 1.3546608450 | -0.5149118351 | 0.1889403276 |
| 409 | 0.3272542486 | 1.4430491930 | 0.1336009891 | 0.1889403276 |
| 410 | -0.3272542486 | 1.4430491930 | 0.1336009891 | 0.1889403276 |
| 411 | 0.5295084972 | 0. | 0.5344039562 | -1.2950162700 |
| 412 | 0.5295084972 | 0. | -0.5344039562 | 1.2950162700 |
| 413 | 0.5295084972 | 0. | 1.3990877210 | -0.0721687833 |
| 414 | -0.5295084972 | 0. | 1.3990877210 | -0.0721687833 |
| 415 | 0.5295084972 | 0.1430153508 | 0.5149118351 | 1.2950162700 |
| 416 | -0.5295084972 | 0.1430153508 | 0.5149118351 | 1.2950162700 |
| 417 | 0.5295084972 | 0.1430153508 | 1.0493157910 | 0.9171356145 |
| 418 | -0.5295084972 | 0.1430153508 | 1.0493157910 | 0.9171356145 |
| 419 | 0.5295084972 | 0.1430153508 | -1.2144556950 | 0.6835925263 |
| 420 | 0.5295084972 | -0.1430153508 | 1.2144556950 | -0.6835925263 |
| 421 | 0.5295084972 | 0.1430153508 | 1.0493157910 | -0.9171356145 |
| 422 | 0.5295084972 | -0.1430153508 | -1.0493157910 | 0.9171356145 |
| 423 | 0.5295084972 | 0.3744190490 | -0.3813108453 | 1.2950162700 |
| 424 | 0.5295084972 | -0.3744190490 | 0.3813108453 | -1.2950162700 |
| 425 | 0.5295084972 | 0.3744190490 | 1.3480566840 | 0.0721687833 |
| 426 | -0.5295084972 | 0.3744190490 | 1.3480566840 | 0.0721687833 |
| 427 | 0.5295084972 | 0.4628073964 | 0.2672019781 | 1.2950162700 |
| 428 | -0.5295084972 | 0.4628073964 | 0.2672019781 | 1.2950162700 |
| 429 | -0.5295084972 | 0.4628073964 | -0.2672019781 | -1.2950162700 |
| 430 | -0.5295084972 | -0.4628073964 | 0.2672019781 | 1.2950162700 |
| 431 | 0.5295084972 | 0.4628073964 | 1.1318857430 | 0.6835925263 |
| 432 | 0.5295084972 | 0.4628073964 | -1.1318857430 | -0.6835925263 |
| 433 | 0.5295084972 | -0.4628073964 | 1.1318857430 | 0.6835925263 |
| 434 | -0.5295084972 | 0.4628073964 | 1.1318857430 | 0.6835925263 |
| 435 | 0.5295084972 | 0.5174343998 | 0.1336009891 | -1.2950162700 |
| 436 | 0.5295084972 | -0.5174343998 | -0.1336009891 | 1.2950162700 |
| 437 | 0.5295084972 | 0.5174343998 | -1.2654867320 | -0.3057118717 |
| 438 | -0.5295084972 | 0.5174343998 | -1.2654867320 | -0.3057118717 |
| 439 | 0.5295084972 | 0.7488380980 | 0.9667458392 | -0.6835925263 |
| 440 | 0.5295084972 | 0.7488380980 | -0.9667458392 | 0.6835925263 |
| 441 | -0.5295084972 | 0.7488380980 | 0.9667458392 | -0.6835925263 |
| 442 | -0.5295084972 | 0.7488380980 | -0.9667458392 | 0.6835925263 |
| 443 | 0.5295084972 | 0.8372264454 | -0.6485128238 | 0.9171356154 |
| 444 | -0.5295084972 | 0.8372264454 | -0.6485128238 | 0.9171356154 |
| 445 | -0.5295084972 | 0.8372264454 | -0.6485128238 | -0.9171356154 |
| 446 | 0.5295084972 | 0.8372264454 | -0.6485128238 | -0.9171356154 |
| 447 | 0.5295084972 | 0.8372264454 | 1.0808547060 | -0.3057118717 |
| 448 | -0.5295084972 | 0.8372264454 | 1.0808547060 | -0.3057118717 |
| 449 | 0.5295084972 | 0.9802417962 | 0.4008029672 | -0.9171356154 |
| 450 | 0.5295084972 | 0.9802417962 | 0.4008029672 | 0.9171356154 |
| 451 | 0.5295084972 | 0.9802417962 | 0.7310827761 | 0.6835925263 |
| 452 | 0.5295084972 | -0.9802417962 | -0.7310827761 | -0.6835925263 |
| 453 | 0.5295084972 | -0.9802417962 | -0.4008029672 | -0.9171356154 |
| 454 | -0.5295084972 | 0.9802417962 | 0.4008029672 | -0.9171356154 |
| 455 | 0.5295084972 | 0.9802417962 | -0.9982847540 | 0.0721687833 |
| 456 | -0.5295084972 | 0.9802417962 | -0.9982847540 | 0.0721687833 |
| 457 | 0.5295084972 | 1.1232571470 | -0.4833729193 | -0.6835925263 |
| 458 | -0.5295084972 | 1.1232571470 | -0.4833729193 | -0.6835925263 |
| 459 | 0.5295084972 | 1.2116454950 | -0.1651399043 | 0.6835925263 |
| 460 | -0.5295084972 | 1.2116454950 | -0.1651399043 | 0.6835925263 |
| 461 | 0.5295084972 | 1.2116454950 | 0.6995438612 | 0.0721687833 |
| 462 | 0.5295084972 | 1.2116454950 | -0.6995438612 | -0.0721687833 |
| 463 | 0.5295084972 | 1.2116454950 | 0.1651399043 | -0.6835925263 |
| 464 | -0.5295084972 | 1.2116454950 | 0.1651399043 | -0.6835925263 |
| 465 | -0.5295084972 | 1.2116454950 | -0.6995438612 | -0.0721687833 |
| 466 | -0.5295084972 | 1.2116454950 | 0.6995438612 | 0.0721687833 |
| 467 | 0.5295084972 | 1.3546608450 | -0.1846320260 | 0.3057118717 |
| 468 | -0.5295084972 | 1.3546608450 | -0.1846320260 | 0.3057118717 |
| 469 | 0.5295084972 | 1.3546608450 | 0.3497719302 | -0.0721687833 |
| 470 | -0.5295084972 | 1.3546608450 | 0.3497719302 | -0.0721687833 |
| 471 | 0.8567627458 | 0. | 0.5344039562 | 1.1060759430 |
| 472 | 0.8567627458 | 0. | -0.5344039562 | -1.1060759430 |
| 473 | 0.8567627458 | 0. | 0.8646837645 | 0.8725328536 |
| 474 | 0.8567627458 | 0. | -0.8646837645 | -0.8725328536 |
| 475 | 0.8567627458 | 0.1430153508 | 0.5149118351 | -1.1060759430 |
| 476 | 0.8567627458 | -0.1430153508 | -0.5149118351 | 1.1060759430 |
| 477 | 0.8567627458 | 0.1430153508 | -1.2144556950 | 0.1167715443 |

| 478 | 0.8567627458 | -0.1430153508 | 1.2144556950 | -0.1167715443 |
| --- | --- | --- | --- | --- |
| 479 | 0.8567627458 | 0.2314036982 | -1.1003468280 | 0.4946521989 |
| 480 | 0.8567627458 | -0.2314036982 | 1.1003468280 | -0.4946521989 |
| 481 | 0.8567627458 | 0.2314036982 | 0.8331448501 | -0.8725328536 |
| 482 | -0.8567627458 | 0.2314036982 | 0.8331448501 | -0.8725328536 |
| 483 | 0.8567627458 | 0.3744190490 | -0.9157148022 | -0.7281952873 |
| 484 | 0.8567627458 | -0.3744190490 | 0.9157148022 | 0.7281952873 |
| 485 | 0.8567627458 | 0.3744190490 | -0.3813108453 | -1.1060759430 |
| 486 | -0.8567627458 | 0.3744190490 | -0.3813108453 | -1.1060759430 |
| 487 | 0.8567627458 | 0.4628073964 | -0.2672019781 | 1.1060759430 |
| 488 | 0.8567627458 | 0.4628073964 | 0.2672019781 | -1.1060759430 |
| 489 | 0.8567627458 | -0.4628073964 | -0.2672019781 | 1.1060759430 |
| 490 | -0.8567627458 | 0.4628073964 | -0.2672019781 | 1.1060759430 |
| 491 | 0.8567627458 | 0.4628073964 | 1.1318857430 | 0.1167715443 |
| 492 | 0.8567627458 | -0.4628073964 | 1.1318857430 | 0.1167715443 |
| 493 | 0.8567627458 | 0.4628073964 | -1.1318857430 | -0.1167715443 |
| 494 | -0.8567627458 | 0.4628073964 | 1.1318857430 | 0.1167715443 |
| 495 | 0.8567627458 | 0.5174343998 | 0.1336009891 | 1.1060759430 |
| 496 | 0.8567627458 | -0.5174343998 | -0.1336009891 | -1.1060759430 |
| 497 | 0.8567627458 | 0.5174343998 | 0.9982847540 | 0.4946521989 |
| 498 | 0.8567627458 | -0.5174343998 | -0.9982847540 | -0.4946521989 |
| 499 | 0.8567627458 | 0.6058227472 | 0.7821138131 | -0.7281952873 |
| 500 | -0.8567627458 | 0.6058227472 | 0.7821138131 | -0.7281952873 |
| 501 | 0.8567627458 | 0.6058227472 | -0.6169739091 | -0.8725328536 |
| 502 | -0.8567627458 | 0.6058227472 | -0.6169739091 | -0.8725328536 |
| 503 | 0.8567627458 | -0.6058227472 | 0.9472537170 | -0.4946521989 |
| 504 | 0.8567627458 | 0.6058227472 | -0.9472537170 | 0.4946521989 |
| 505 | 0.8567627458 | 0.7488380980 | 0.4323418823 | -0.8725328536 |
| 506 | 0.8567627458 | 0.7488380980 | -0.4323418823 | 0.8725328536 |
| 507 | -0.8567627458 | 0.7488380980 | 0.4323418823 | -0.8725328536 |
| 508 | -0.8567627458 | 0.7488380980 | -0.4323418823 | 0.8725328536 |
| 509 | 0.8567627458 | 0.7488380980 | 0.9667458392 | -0.1167715443 |

| 510 | 0.8567627458 | 0.7488380980 | -0.9667458392 | 0.11677154430 |
|---|---|---|---|---|
| 511 | -0.8567627458 | 0.7488380980 | 0.9667458392 | -0.1167715443 |
| 512 | -0.8567627458 | 0.7488380980 | -0.9667458392 | 0.1167715443 |
| 513 | 0.8567627458 | 0.8372264454 | 0.2161709413 | 0.8725328536 |
| 514 | -0.8567627458 | 0.8372264454 | 0.2161709413 | 0.8725328536 |
| 515 | 0.8567627458 | 0.8372264454 | 0.7505748982 | 0.4946521989 |
| 516 | -0.8567627458 | 0.8372264454 | 0.7505748982 | 0.4946521989 |
| 517 | 0.8567627458 | 0.9802417962 | -0.1336009891 | 0.7281952873 |
| 518 | -0.8567627458 | 0.9802417962 | -0.1336009891 | 0.7281952873 |
| 519 | 0.8567627458 | 0.9802417962 | 0.7310827761 | 0.1167715443 |
| 520 | -0.8567627458 | 0.9802417962 | 0.7310827761 | 0.1167715443 |
| 521 | 0.8567627458 | 1.0686301440 | -0.3497719302 | -0.4946521989 |
| 522 | -0.8567627458 | 1.0686301440 | -0.3497719302 | -0.4946521989 |
| 523 | 0.8567627458 | 1.1232571470 | -0.4833729193 | -0.1167715443 |
| 524 | -0.8567627458 | 1.1232571470 | -0.4833729193 | -0.1167715443 |
| 525 | 0.8567627458 | 1.1232571470 | 0.0510310370 | -0.4946521989 |
| 526 | -0.8567627458 | 1.1232571470 | 0.0510310370 | -0.4946521989 |
| 527 | 0.8567627485 | 1.2116454950 | 0.1651399043 | -0.1167715443 |
| 528 | 0.8567627458 | 1.2116454950 | -0.1651399043 | 0.1167715443 |
| 529 | 0.8567627458 | -1.2116454950 | -0.1651399043 | 0.1167715443 |
| 530 | -0.8567627458 | 1.2116454950 | -0.1651399043 | 0.1167715443 |
| 531 | 1.059019944 | 0. | 0.8646837645 | -0.6114237430 |
| 532 | 1.0590169941 | 0. | -0.8646837645 | 0.6114237430 |
| 533 | 1.0590169944 | 0.1430153508 | -0.3497719302 | 0.9893043979 |
| 534 | 1.0590169944 | -0.1430153508 | 0.3497719302 | -0.9893043979 |
| 535 | 1.0590169944 | 0.1430153508 | 1.0493157910 | 0. |
| 536 | 1.0590169944 | -0.1430153508 | -1.0493157910 | 0. |
| 537 | 1.0590169944 | 0.2314036982 | 0.2987408932 | 0.9893043979 |
| 538 | 1.0590169944 | -0.2314036982 | -0.2987408932 | -0.9893043979 |
| 539 | 1.0590169944 | 0.2314036982 | 0.8331448501 | 0.6114237430 |
| 540 | 1.0590169944 | -0.2314036982 | -0.8331448501 | -0.6114237430 |
| 541 | 1.0590169944 | 0.3744190490 | -0.0510310370 | -0.9893043979 |

| 542 | 1.0590169944 | -0.3744190490 | 0.0510310370 | 0.9893043979 |
|---|---|---|---|---|
| 543 | 1.0590169944 | 0.3744190490 | -0.9157148022 | -0.3778806546 |
| 544 | 1.0590169944 | -0.3744190490 | 0.9157148022 | 0.3778806546 |
| 545 | 1.0590169944 | 0.6058227472 | -0.6169739091 | 0.6114237430 |
| 546 | -1.0590169944 | 0.6058227472 | -0.6169739091 | 0.6114237430 |
| 547 | 1.0590169944 | 0.6058227472 | 0.7821138131 | -0.3778806546 |
| 548 | -1.0590169944 | 0.6058227472 | 0.7821138131 | -0.3778806546 |
| 549 | 1.0590169944 | 0.7488380980 | 0.4323418823 | 0.6114237430 |
| 550 | 1.0590169944 | 0.7488380980 | -0.4323418823 | -0.6114237430 |
| 551 | 1.0590169944 | -0.7488380980 | -0.4323418823 | -0.6114237430 |
| 552 | 1.0590169944 | -0.7488380980 | -0.4323418823 | -0.6114237430 |
| 553 | 1.0590169944 | 0.8372264454 | 0.2161709413 | -0.6114237430 |
| 554 | -1.0590169944 | 0.8372264454 | 0.2161709413 | -0.6114237430 |
| 555 | 1.0590169944 | 0.8372264454 | -0.6485128238 | 0. |
| 556 | -1.0590169944 | 0.8372264454 | -0.6485128238 | 0. |
| 557 | 1.0590169944 | 0.9802417962 | -0.1336009891 | 0.3778806546 |
| 558 | -1.0590169944 | 0.9802417962 | -0.1336009891 | 0.3778806546 |
| 559 | 1.0590169944 | 0.9802417962 | 0.4008029672 | 0. |
| 560 | -1.0590169944 | 0.9802417962 | 0.4008029672 | 0. |
| 561 | 1.1840169944 | 0. | 0. | 0.9171356154 |
| 562 | 1.1840169944 | 0. | 0. | -0.9171356154 |
| 563 | 1.1840169944 | 0. | 0.8646837645 | 0.3057118717 |
| 564 | 1.1840169944 | 0. | -0.8646837645 | -0.3057118717 |
| 565 | 1.1840169944 | 0.2314036982 | 0.8331448501 | -0.3057118717 |
| 566 | 1.1840169944 | -0.2314036982 | -0.8331448501 | 0.3057118717 |
| 567 | 1.1840169944 | 0.2314036982 | -0.5659428717 | 0.6835925263 |
| 568 | 1.1840169944 | -0.2314036982 | 0.5659428717 | -0.6835925263 |
| 569 | 1.1840169944 | 0.3744190490 | 0.4833729193 | 0.6835925263 |
| 570 | 1.1840169944 | -0.3744190490 | -0.4833729193 | -0.6835925263 |
| 571 | 1.1840169944 | 0.6058227472 | -0.0825688521 | -0.6835925263 |
| 572 | 1.1840169944 | -0.6058227472 | 0.0825688521 | 0.6835925263 |
| 573 | 1.1840169944 | 0.6058227472 | -0.6169739091 | -0.3057118717 |
| 574 | -1.1840169944 | 0.6058227472 | -0.6169739091 | -0.3057118717 |
| 575 | 1.1840169944 | 0.7488380980 | -0.4323418823 | 0.3057118717 |
| 576 | 1.1840169944 | 0.7488380980 | 0.4323418823 | -0.3057118717 |
| 577 | 1.1840169944 | -0.7488380980 | -0.4323418823 | 0.3057118717 |
| 578 | 1.1840169944 | 0.7488380980 | -0.4323418823 | 0.3057118717 |
| 579 | 1.1840169944 | 0.8372264454 | 0.2161709413 | 0.3057118717 |
| 580 | 1.1840169944 | -0.8372264454 | -0.2161709413 | -0.3057118717 |
| 581 | 1.3862712430 | 0. | 0. | 0.5668209828 |
| 582 | 1.3862712430 | 0. | 0. | -0.5668209828 |
| 583 | 1.3862712430 | 0. | 0.5344039562 | 0.1889403276 |
| 584 | 1.3862712430 | 0. | -0.5344039562 | -0.1889403276 |
| 585 | 1.3862712430 | 0.1430153508 | -0.3497719302 | 0.4224834156 |
| 586 | 1.3862712430 | -0.1430153508 | 0.3497719302 | -0.4224834156 |
| 587 | 1.3862712430 | 0.1430153508 | 0.5149118351 | -0.1889403276 |
| 588 | 1.3862712430 | -0.1430153508 | -0.5149118351 | 0.1889403276 |
| 589 | 1.3862712430 | 0.2314036982 | 0.2987408932 | 0.4224834156 |
| 590 | 1.3862712430 | -0.2314036982 | -0.2987408932 | -0.4224834156 |
| 591 | 1.3862712430 | 0.3744190490 | -0.0510310370 | -0.4224834156 |
| 592 | 1.3862712430 | -0.3744190490 | 0.0510310370 | 0.4224834156 |
| 593 | 1.3862712430 | 0.3744190490 | -0.3813108453 | -0.1889403276 |
| 594 | 1.3862712430 | -0.3744190490 | 0.3813108453 | 0.1889403276 |
| 595 | 1.3862712430 | 0.4628073964 | 0.2672019781 | -0.1889403276 |
| 596 | 1.3862712430 | 0.4628073964 | -0.2672019781 | 0.1889403276 |
| 597 | 1.3862712430 | -0.4628073964 | 0.2672019781 | -0.1889403276 |
| 598 | 1.3862712430 | -0.4628073964 | -0.2672019871 | 0.1889403276 |
| 599 | 1.3862712430 | 0.5174343998 | 0.1336009891 | 0.1889403276 |
| 600 | 1.3862712430 | -0.5174343998 | -0.1336009891 | -0.1889403276 |

# The Joint Relations between the Vertices of the 120-Cell

| | | |
|---|---|---|
| 1(61,62,302,304) | 31(32,56,57,65) | 61(1,39,115,117) |
| 2(11,301,357,358) | 32(31,58,66,355) | 62(1,40,116,118) |
| 3(6,8,51,52) | 33(34,60,63,67) | 63(7,33,113,438) |
| 4(12,59,301,360) | 34(33,64,68,359) | 64(7,34,114,437) |
| 5(7,53,309,354) | 35(36,72,73,352) | 65(10,31,123,127) |
| 6(3,16,55,356) | 36(35,70,74,351) | 66(10,32,128,424) |
| 7(5,63,64,315) | 37(53,71,75,338) | 67(13,33,125,131) |
| 8(3,19,71,369) | 38(54,69,337,376) | 68(13,34,126,134) |
| 9(21,70,72,305) | 39(40,61,77,81) | 69(38,122,308,441) |
| 10(17,18,65,66) | 40(39,62,79,83) | 70(9,36,142,420) |
| 11(2,12,78,80) | 41(42,78,88,102) | 71(8,37,121,139) |
| 12(4,11,82,84) | 42(41,80,87,104) | 72(9,35,119,140) |
| 13(14,20,67,68) | 43(44,82,85,98) | 73(17,35,123,143) |
| 14(13,23,75,76) | 44(43,84,86,100) | 74(17,36,144,424) |
| 15(21,85,86,307) | 45(48,91,95,105) | 75(14,37,125,147) |
| 16(6,19,87,88) | 46(47,89,93,107) | 76(14,126,148,338) |
| 17(10,22,73,74) | 47(46,90,94,108) | 77(18,39,127,150) |
| 18(10,24,77,79) | 48(45,92,96,106) | 78(11,41,146,430) |
| 19(8,16,91,92) | 49(50,97,101,109) | 79(18,40,128,453) |
| 20(13,25,81,83) | 50(49,99,103,110) | 80(11,42,129,145) |
| 21(9,15,89,90) | 51(3,111,121,336) | 81(20,39,131,151) |
| 22(17,24,93,94) | 52(3,335,412,422) | 82(12,43,132,146) |
| 23(14,25,95,96) | 53(5,37,113,120) | 83(20,40,134,452) |
| 24(18,22,101,103) | 54(38,119,305,414) | 84(12,44,145,433) |
| 25(20,23,97,99) | 55(6,111,135,332) | 85(15,43,137,155) |
| 26(27,29,98,100) | 56(31,112,136,306) | 86(15,44,138,156) |
| 27(26,28,102,104) | 57(31,115,302,429) | 87(16,42,154,436) |
| 28(27,30,105,106) | 58(32,116,130,302) | 88(16,41,135,149) |
| 29(26,30,107,108)) | 59(4,132,334,418) | 89(21,46,140,155) |
| 30(28,29,109,110) | 60(33,117,133,304) | 90(21,47,142,156) |
| 91(19,45,139,149) | 123(65,73,112,187) | 155(85,89,162,210) |
| 92(19,48,141,154) | 124(111,366,374,490) | 156(86,90,165,212) |
| 93(22,46,143,159) | 125(67,75,113,191) | 157(98,102,146,221) |
| 94(22,47,144,160) | 126(68,76,114,194) | 158(100,104,145,222) |
| 95(23,45,147,161) | 127(65,77,115,195) | 159(93,101,167,217) |
| 96(23,48,148,166) | 128(66,79,116,496) | 160(94,103,168,218) |
| 97(25,49,151,161) | 129(80,186,357,436) | 161(95,97,169,219) |
| 98(26,43,157,162) | 130(58,378,435,485) | 162(98,107,155,223) |
| 99(25,50,166,452) | 131(67,81,117,197) | 163(102,105,149,225) |
| 100(26,44,158,165) | 132(59,82,137,183) | 164(104,106,154,226) |
| 101(24,49,150,159) | 133(60,184,384,438) | 165(100,108,156,224) |
| 102(27,41,157,163) | 134(68,83,118,498) | 166(96,99,170,220) |
| 103(24,50,160,453) | 135(55,88,188,430) | 167(107,109,159,228) |
| 104(27,42,158, 164) | 136(56,189,387,429) | 168(108,110,160,230) |
| 105(28,45,163,169) | 137(85,132,193,364) | 169(105,109,161,227) |
| 106(28,48,164,170) | 138(86,363,433,492) | 170(106,110,166,529) |
| 107(29,46,162,167) | 139(71,91,147,199) | 171(115,173,237,486) |
| 108(29,47,165,168) | 140(72,89,143,204) | 172(174,185,238,416) |
| 109(30,49,167,169) | 141(92,148,200,369) | 173(117,171,184,239) |
| 110(30,50,168,170) | 142(70,90,144,503) | 174(172,183,240,418) |
| 111(51,55,124,175) | 143(73,93,140,206) | 175(111,181,188,234) |
| 112(56,123,176,352) | 144(74,94,142,208) | 176(112,189,233,482) |
| 113(53,63,125,178) | 145(80,84,158,202) | 177(179,193,236,414) |
| 114(64,126,354,477) | 146(78,82,157,201) | 178(113,180,192,235) |
| 115(57,61,127,171) | 147(75,95,139,209) | 179(119,177,204,232) |
| 116(58,62,128,472) | 148(76,96,141,211) | 180(120,178,203,231) |
| 117(60,61,131,173) | 149(88,91,163,205) | 181(121,175,199,231) |
| 118(62,134,359,474) | 150(77,101,151,213) | 182(200,422,476,532) |
| 119(54,72,122,179) | 151(81,97,150,215) | 183(132,174,201,243) |
| 120(53,121,180,370) | 152(153,383,399,516) | 184(133,173,244,502) |
| 121(51,71,120,181) | 153(152,379,403,514) | 185(172,201,241,430) |
| 122(69,119,352,482) | 154(87,92,164,207) | 186(129,202,471,542) |
| 187(123,195,206,233) | 219(161,209,215,259) | 251(270,280,514,516) |
| 188(135,175,205,241) | 220(166,211,216,260) | 252(272,502,522,574) |
| 189(136,176,242,507) | 221(157,223,225,250) | 253(205,225,271,276) |
| 190(208,424,496,534) | 222(224,226,158,552) | 254(207,226,572,577) |
| 191(125,197,209,235) | 223(162,221,228,255) | 255(210,223,273,275) |
| 192(178,244,438,512) | 224(222,165,230,256) | 256(212,224,274,278) |
| 193(137,177,210,243) | 225(163,221,227,253) | 257(217,228,275,279) |
| 194(126,211,498,536) | 226(164,222,254,529) | 258(218,230,278,580) |
| 195(127,187,213,237) | 227(225,228,169,259) | 259(219,227,276,279) |
| 196(238,428,490,514) | 228(167,223,227,257) | 260(220,529,577,580) |
| 197(131,191,215,239) | 229(470,526,530,560) | 261(233,237,242,281) |
| 198(240,434,494,516) | 230(168,224,258,529) | 262(234,238,241,282) |
| 199(139,181,205,247) | 231(180,181,265,268) | 263(235,239,244,283) |
| 200(141,182,207,248) | 232(179,266,267,482) | 264(236,240,243,284) |
| 201(146,183,185,250) | 233(187,176,261,267) | 265(231,235,247,287) |
| 202(145,186,484,552) | 234(175,262,268,490) | 266(232,236,288,548) |
| 203(180,442,512,546) | 235(178,191,263,265) | 267(232,233,245,285) |
| 204(140,179,210,245) | 236(177,264,266,494) | 268(231,234,286,546) |
| 205(149,188,199,253) | 237(171,195,261,269) | 269(237,239,249,289) |
| 206(143,187,217,245) | 238(172,196,262,270) | 270(238,240,251,290) |
| 207(154,200,254,489) | 239(173,197,263,269) | 271(241,250,253,291) |
| 208(144,190,218,246) | 240(174,198,264,270) | 272(242,252,292,554) |
| 209(147,191,219,247) | 241(185,188,262,271) | 273(243,250,255,293) |
| 210(155,193,204,255) | 242(189,261,272,486) | 274(256,544,552,594) |
| 211(148,194,220,248) | 243(183,193,264,273) | 275(245,255,257,296) |
| 212(156,256,492,503) | 244(184,192,263,574) | 276(247,253,259,295) |
| 213(150,195,217,249) | 245(204,206,267,275) | 277(298,548,554,560) |
| 214(218,453,496,551) | 246(208,278,503,568) | 278(246,256,258,597) |
| 215(151,197,219,249) | 247(199,209,265,276) | 279(249,257,259,299) |
| 216(220,452,498,551) | 248(200,211,566,577) | 280(251,300,558,560) |
| 217(159,206,213,257) | 249(213,215,269,279) | 281(261,285,289,292) |
| 218(160,208,214,258) | 250(201,221,271,273) | 282(262,286,290,291) |
| 283(263,287,289,294) | 315(7,321,385,386) | 347(346,390,394,408) |
| 284(264,288,290,293) | 316(306,319,387,388) | 348(345,392,396,406) |
| 285(267,281,288,296) | 317(310,322,373,374) | 349(350,397,401,409) |
| 286(268,282,287,297) | 318(310,324,377,379) | 350(349,399,403,410) |
| 287(265,283,286,295) | 319(308,316,391,392) | 351(36,303,411,421) |
| 288(266,284,285,298) | 320(313,325,381,383) | 352(35,112,122,303) |
| 289(269,281,283,299) | 321(309,315,389,390) | 353(305,337,413,420) |
| 290(270,282,284,300) | 322(317,324,393,394) | 354(5,114,338,419) |
| 291(271,282,293,295) | 323(314,325,395,396) | 355(32,306,411,435) |
| 292(272,281,294,298) | 324(318,322,401,403) | 356(6,331,412,436) |
| 293(273,284,291,296) | 325(320,323,397,399) | 357(2,129,331,415) |
| 294(283,292,297,574) | 326(327,329,398,400) | 358(2,332,416,430) |
| 295(276,287,291,299) | 327(326,328,402,404) | 359(34,118,304,432) |
| 296(275,285,293,299) | 328(327,330,405,406) | 360(4,333,417,433) |
| 297(286,578,294,300) | 329(326,330,407,408)) | 361(301,339,415,417) |
| 298(277,288,292,300) | 330(328,329,409,410) | 362(301,340,416,418) |
| 299(279,289,295,296) | 331(332,356,357,365) | 363(138,307,333,413) |
| 300(280,290,297,298) | 332(55,331,358,366) | 364(137,307,334,414) |
| 301(2,4,361,362) | 333(334,360,363,367) | 365(310,331,423,427) |
| 302(1,57,58,311) | 334(59,333,364,368) | 366(124,310,332,428) |
| 303(306,308,351,352) | 335(52,336,372,373) | 367(313,333,425,431) |
| 304(1,60,312,359) | 336(51,335,370,374) | 368(313,334,426,434) |
| 305(9,54,307,353) | 337(38,353,371,375) | 369(8,141,338,422) |
| 306(56,303,316,355) | 338(37,76,354,369) | 370(120,309,336,442) |
| 307(15,305,363,364) | 339(340,361,377,381) | 371(308,337,421,439) |
| 308(69,303,319,371) | 340(339,362,379,383) | 372(309,335,419,440) |
| 309(5,321,370,372) | 341(342,378,388,402) | 373(317,335,423,443) |
| 310(317,318,365,366) | 342(341,380,387,404) | 374(124,317,336,444) |
| 311(302,312,378,380) | 343(344,382,385,398) | 375(314,337,425,447) |
| 312(304,311,382,384) | 344(343,384,386,400) | 376(38,314,426,448) |
| 313(314,320,367,368) | 345(348,391,395,405) | 377(318,339,427,450) |
| 314(313,323,375,376) | 346(347,389,393,407) | 378(130,311,341,446) |
| **379(153,318,340,428)** | **411(351,355,424,475)** | **443(373,393,440,506)** |
| **380(311,342,429,445)** | **412(52,356,423,476)** | **444(374,394,442,508)** |
| **381(320,339,431,451)** | **413(353,363,425,478)** | **445(380,384,458,502)** |
| **382(312,343,432,446)** | **414(54,177,364,426)** | **446(378,382,457,501)** |
| **383(152,320,340,434)** | **415(357,361,427,471)** | **447(375,395,439,509)** |
| **384(133,312,344,445)** | **416(172,358,362,428)** | **448(376,396,441,511)** |
| **385(315,343,437,455)** | **417(360,361,431,473)** | **449(388,391,463,505)** |
| **386(315,344,438,456)** | **418(59,174,362,434)** | **450(377,401,451,513)** |
| **387(136,316,342,454)** | **419(354,372,422,479)** | **451(381,397,450,515)** |
| **388(316,341,435,449)** | **420(70,353,421,480)** | **452(83,99,216,453)** |
| **389(321,346,440,455)** | **421(351,371,420,481)** | **453(79,103,214,452)** |
| **390(321,347,442,456)** | **422(52,182,369,419)** | **454(387,392,464,507)** |
| **391(319,345,439,449)** | **423(365,373,412,487)** | **455(385,389,462,510)** |
| **392(319,348,441,454)** | **424(66,74,190,411)** | **456(386,390,465,512)** |
| **393(322,346,443,459)** | **425(367,375,413,491)** | **457(398,402,446,521)** |
| **394(322,347,444,460)** | **426(368,376,414,494)** | **458(400,404,445,522)** |
| **395(323,345,447,461)** | **427(365,377,415,495)** | **459(393,401,467,517)** |
| **396(323,348,448,466)** | **428(196,366,379,416)** | **460(394,403,468,518)** |
| **397(325,349,451,461)** | **429(57,136,380,486)** | **461(395,397,469,519)** |
| **398(326,343,457,462)** | **430(78,135,185,358)** | **462(398,407,455,523)** |
| **399(152,325,350,466)** | **431(367,381,417,497)** | **463(402,405,449,525)** |
| **400(326,344,458,465)** | **432(359,382,437,483)** | **464(404,406,454,526)** |
| **401(324,349,450,459)** | **433(84,138,360,484)** | **465(400,408,456,524)** |
| **402(327,341,457,463)** | **434(198,368,383,418)** | **466(396,399,470,520)** |
| **403(153,324,350,460)** | **435(130,355,388,488)** | **467(407,409,459,528)** |
| **404(327,342,458,464)** | **436(87,129,356,489)** | **468(408,410,460,530)** |
| **405(328,345,463,469)** | **437(64,385,432,493)** | **469(405,409,461,527)** |
| **406(328,348,464,470)** | **438(63,133,192,386)** | **470(229,406,410,466)** |
| **407(329,346,462,467)** | **439(371,391,447,499)** | **471(186,415,473,537)** |
| **408(329,347,465,468)** | **440(372,389,443,504)** | **472(116,474,485,538)** |
| **409(330,349,467,469)** | **441(69,392,448,500)** | **473(417,471,484,539)** |
| **410(330,350,468,470)** | **442(203,370,390,444)** | **474(118,472,483,540)** |
| 475(411,481,488,534) | 507(189,454,500,554) | 539(473,497,563,569) |
| 476(182,412,489,533) | 508(444,490,518,546) | 540(474,498,564,570) |
| 477(114,479,493,536) | 509(447,491,519,547) | 541(485,488,562,571) |
| 478(413,480,492,535) | 510(455,493,504,555) | 542(186,489,561,572) |
| 479(419,477,504,532) | 511(448,494,520,548) | 543(483,493,564,573) |
| 480(420,478,503,531) | 512(192,203,456,556) | 544(274,484,492,563) |
| 481(421,475,499,531) | 513(450,495,517,549) | 545(504,506,567,575) |
| 482(122,176,232,500) | 514(153,196,251,518) | 546(203,268,508,578) |
| 483(432,474,501,543) | 515(451,497,519,549) | 547(499,509,565,576) |
| 484(202,433,473,544) | 516(152,198,251,520) | 548(266,277,500,511) |
| 485(130,472,501,541) | 517(459,506,513,557) | 549(513,515,569,579) |
| 486(171,242,429,502) | 518(460,508,514,558) | 550(501,521,571,573) |
| 487(423,495,506,533) | 519(461,509,515,559) | 551(214,216,570,580) |
| 488(435,475,505,541) | 520(466,511,516,560) | 552(202,222,274,572) |
| 489(207,436,476,542) | 521(457,523,525,550) | 553(505,525,571,576) |
| 490(124,196,234,508) | 522(252,524,526,458) | 554(272,277,507,526) |
| 491(425,497,509,535) | 523(462,521,528,555) | 555(510,523,573,575) |
| 492(138,212,478,544) | 524(522,465,530,556) | 556(512,524,574,578) |
| 493(437,477,510,543) | 525(463,521,527,553) | 557(517,528,575,579) |
| 494(198,236,426,511) | 526(229,464,522,554) | 558(280,518,530,578) |
| 495(427,487,513,537) | 527(469,525,528,559) | 559(519,527,576,579) |
| 496(128,190,214,538) | 528(467,523,527,557) | 560(229,277,280,520) |
| 497(431,491,515,539) | 529(170,226,230,260) | 561(533,537,542,581) |
| 498(134,194,216,540) | 530(229,468,524,558) | 562(534,538,541,582) |
| 499(439,481,505,547) | 531(480,481,565,568) | 563(535,539,544,583) |
| 500(441,482,507,548) | 532(182,479,566,567) | 564(536,540,543,584) |
| 501(446,483,485,550) | 533(476,487,561,567) | 565(531,535,547,587) |
| 502(184,252,445,486) | 534(190,475,562,568) | 566(248,532,536,588) |
| 503(142,212,246,480) | 535(478,491,563,565) | 567(532,533,545,585) |
| 504(440,479,510,545) | 536(194,477,564,566) | 568(246,531,534,586) |
| 505(449,488,499,553) | 537(471,495,561,569) | 569(537,539,549,589) |
| 506(443,487,517,545) | 538(472,496,562,570) | 570(538,540,551,590) |
| 571(541,550,553,591) | 581(561,585,589,592) | 591(571,582,593,595) |
| 572(254,542,552,592) | 582(562,586,590,591) | 592(572,581,594,598) |
| 573(543,550,555,593) | 583(563,587,589,594) | 593(573,584,591,596) |
| 574(244,252,294,556) | 584(564,588,590,593) | 594(274,583,592,597) |
| 575(545,555,557,596) | 585(567,581,588,596) | 595(576,587,591,599) |
| 576(547,553,559,595) | 586(568,582,587,597) | 596(575,585,593,599) |
| 577(248,254,260,598) | 587(565,583,586,595) | 597(586,278,594,600) |
| 578(297,546,556,558) | 588(566,584,585,598) | 598(577,588,592,600) |
| 579(549,557,559,599) | 589(569,581,583,599) | 599(579,589,595,596) |
| 580(258,260,551,600) | 590(570,582,584,600) | 600(580,590,597,598) |

# The 4-dimensional Coordinates of Vertices of the 600-Cell

| No | $x_1$ | $x_2$ | $x_3$ | $x_4$ |
|----|-------|-------|-------|-------|
| 1  | 0. | 0. | 1.3211192340 | -0.9341723545 |
| 2  | 0. | -0.2185080123 | 0.5344039563 | -1.5115226210 |
| 3  | 0. | 0.2185080123 | 1.6032118690 | 0. |
| 4  | 0. | -0.3535533906 | -1.2729320600 | -0.9341723545 |
| 5  | 0. | -0.3535533906 | -0.4564354698 | -1.5115226210 |
| 6  | 0. | -0.5720614029 | 1.3990877210 | 0.5773502668 |
| 7  | 0. | 0.5720614029 | -0.0779684865 | -1.5115226210 |
| 8  | 0. | -0.9256147931 | 0.9426522510 | -0.9341723545 |
| 9  | 0. | 0.9256147931 | 1.1949635730 | -0.5773502668 |
| 10 | 0. | -1.1441228050 | -0.6605596176 | -0.9341723545 |
| 11 | 0. | -1.1441228050 | 0.6605596176 | 0.9341723545 |
| 12 | 0. | -1.2791681840 | 0.9908394253 | 0. |
| 13 | 0. | -1.2791681840 | -0.3302798088 | 0.9341723545 |
| 14 | 0. | -1.4976761961 | -0.6123724428 | 0. |
| 15 | 0. | -1.4976761961 | 0.2041241477 | -0.5773502668 |
| 16 | -0.5 | 0.3535533906 | -0.8646837645 | -1.2228474870 |
| 17 | -0.5 | -0.5720614029 | 1.3990877210 | -0.2886751334 |
| 18 | -0.5 | 0.5720614029 | 0.7385281032 | -1.2228474870 |
| 19 | -0.5 | -0.9256147931 | -1.1949635730 | -0.2886751334 |
| 20 | -0.5 | -0.9256147931 | 0.1261556613 | -1.2228474870 |
| 21 | -0.5 | -1.4976761961 | 0.2041241477 | 0.2886751334 |
| 22 | 0.5 | 0.3535533906 | -0.8646837645 | -1.2228474870 |
| 23 | 0.5 | -0.5720614029 | 1.3990877210 | -0.2886751334 |
| 24 | 0.5 | 0.5720614029 | 0.7385281032 | -1.2228474870 |
| 25 | 0.5 | -0.9256147931 | -1.1949635730 | -0.2886751334 |
| 26 | 0.5 | -0.9256147931 | 0.1261556613 | -1.2228474870 |
| 27 | 0.5 | -1.4976761961 | 0.2041241477 | 0.2886751334 |
| 28 | -0.8090169944 | 0. | 0. | -1.4012585320 |
| 29 | -0.8090169944 | -0.3535533906 | 0.8646837645 | -1.0444364440 |
| 30 | -0.8090169944 | -0.5720614029 | -0.7385281032 | -1.0444364440 |
| 31 | -0.8090169944 | -0.9256147931 | 0.9426522510 | 0.4670861773 |
| 32 | -0.8090169944 | -1.1441228050 | -0.6605596176 | 0.4670861773 |
| 33 | -0.8090169944 | -1.1441228050 | 0.6605596176 | -0.4670861773 |
| 34 | -0.8090169944 | -1.2791681840 | -0.3302798088 | -0.4670861773 |
| 35 | 0.8090169944 | 0. | 0. | -1.4012585320 |
| 36 | 0.8090169944 | 0. | -1.3211192340 | -0.4670861773 |
| 37 | 0.8090169944 | 0. | 1.3211192340 | 0.4670861773 |
| 38 | 0.8090169944 | -0.3535533906 | 0.8646837645 | -1.0444364440 |
| 39 | 0.8090169944 | -0.3535533906 | -1.2729320600 | 0.4670861773 |
| 40 | 0.8090169944 | 0.3535533906 | 1.2729320600 | -0.4670861773 |
| 41 | 0.8090169944 | -0.5720614029 | -0.7385281032 | -1.0444364440 |
| 42 | 0.8090169944 | -0.9256147931 | 0.9426522510 | 0.4670861773 |
| 43 | 0.8090169944 | -0.9256147931 | 0.1261556613 | 1.0444364440 |
| 44 | 0.8090169944 | 0.9256147931 | -0.1261556613 | -1.0444364440 |
| 45 | 0.8090169944 | -1.1441228050 | -0.6605596176 | 0.4670861773 |
| 46 | 0.8090169944 | -1.1441228050 | 0.6605596176 | -0.4670861773 |
| 47 | 0.8090169944 | -1.2791681840 | -0.3302798088 | -0.4670861773 |
| 48 | -1.3090169944 | -0.3535533906 | 0.8646837645 | -0.1784110442 |
| 49 | -1.3090169944 | -0.5720614029 | 0.0779684865 | -0.7557613103 |
| 50 | 1.3090169944 | -0.2185080123 | 0.5344039563 | 0.7557613103 |
| 51 | 1.3090169944 | 0.2185080123 | -0.5344039563 | -0.7557613103 |
| 52 | 1.3090169944 | -0.3535533906 | -0.4564354698 | 0.7557613103 |
| 53 | 1.3090169944 | -0.3535533906 | 0.8646837645 | -0.1784110442 |
| 54 | 1.3090169944 | 0.3535533906 | 0.4564354698 | -0.7557613103 |
| 55 | 1.3090169944 | -0.5720614029 | -0.7385281032 | -0.1784110442 |
| 56 | 1.3090169944 | -0.5720614029 | 0.0779684865 | -0.7557613103 |
| 57 | 1.3090169944 | 0.5720614029 | 0.7385281032 | 0.1784110442 |
| 58 | 1.3090169944 | -0.9256147931 | 0.1261556613 | 0.1784110442 |
| 59 | 1.3090169944 | 0.9256147931 | -0.1261556613 | -0.1784110442 |
| 60 | 1.6180339891 | 0. | 0. | 0. |
| 61 | 0. | 0. | -1.3211192340 | 0.9341723545 |
| 62 | 0. | 0.2185080123 | -0.5344039563 | 1.5115226210 |
| 63 | 0. | -0.2185080123 | -1.6032118690 | 0. |
| 64 | 0. | 0.3535533906 | 1.2729320600 | 0.9341723545 |
| 65 | 0. | 0.3535533906 | 0.4564354698 | 1.5115226210 |
| 66 | 0. | 0.5720614029 | -1.3990877210 | -0.5773502668 |
| 67 | 0. | -0.5720614029 | 0.0779684865 | 1.5115226210 |
| 68 | 0. | 0.9256147931 | -0.9426522510 | 0.9341723545 |
| 69 | 0. | -0.9256147931 | -1.1949635730 | 0.5773502668 |
| 70 | 0. | 1.1441228050 | 0.6605596176 | 0.9341723545 |
| 71 | 0. | 1.1441228050 | -0.6605596176 | -0.9341723545 |
| 72 | 0. | 1.2791681840 | -0.9908394253 | 0. |
| 73 | 0. | 1.2791681840 | 0.3302798088 | -0.9341723545 |
| 74 | 0. | 1.4976761961 | 0.6123724428 | 0. |
| 75 | 0. | 1.4976761961 | -0.2041241477 | 0.5773502668 |
| 76 | 0.5 | -0.3535533906 | 0.8646837645 | 1.2228474870 |
| 77 | 0.5 | 0.5720614029 | -1.3990877210 | 0.2886751334 |
| 78 | 0.5 | -0.5720614029 | -0.7385281032 | 1.2228474870 |
| 79 | 0.5 | 0.9256147931 | 1.1949635730 | 0.2886751334 |
| 80 | 0.5 | 0.9256147931 | -0.1261556613 | 1.2228474870 |
| 81 | 0.5 | 1.4976761961 | -0.2041241477 | -0.2886751334 |
| 82 | -0.5 | -0.3535533906 | 0.8646837645 | 1.2228474870 |
| 83 | -0.5 | 0.5720614029 | -1.3990877210 | 0.2886751334 |
| 84 | -0.5 | -0.5720614029 | -0.7385281032 | 1.2228474870 |
| 85 | -0.5 | 0.9256147931 | 1.1949635730 | 0.2886751334 |
| 86 | -0.5 | 0.9256147931 | -0.1261556613 | 1.2228474870 |
| 87 | -0.5 | 1.4976761961 | -0.2041241477 | -0.2886751334 |
| 88 | 0.8090169944 | 0. | 0. | 1.4012585320 |
| 89 | 0.8090169944 | 0.3535533906 | -0.8646837645 | 1.0444364440 |
| 90 | 0.8090169944 | 0.5720614029 | 0.7385281032 | 1.0444364440 |
| 91 | 0.8090169944 | 0.9256147931 | -0.9426522510 | -0.4670861773 |
| 92 | 0.8090169944 | 1.1441228050 | 0.6605596176 | -0.4670861773 |
| 93 | 0.8090169944 | 1.1441228050 | -0.6605596176 | 0.4670861773 |

| 94  | 0.8090169944  | 1.2791681840  | 0.3302798088  | 0.4670861773  |
| --- | --- | --- | --- | --- |
| 95  | -0.8090169944 | 0.            | 0.            | 1.4012585320  |
| 96  | -0.8090169944 | 0.            | 1.3211192340  | 0.4670861773  |
| 97  | -0.8090169944 | 0.            | -1.3211192340 | -0.4670861773 |
| 98  | -0.8090169944 | 0.3535533906  | -0.8646837645 | 1.0444364440  |
| 99  | -0.8090169944 | 0.3535533906  | 1.2729320600  | -0.4670861773 |
| 100 | -0.8090169944 | -0.3535533906 | -1.2729320600 | 0.4670861773  |
| 101 | -0.8090169944 | 0.5720614029  | 0.7385281032  | 1.0444354440  |
| 102 | -0.8090169944 | 0.9256147931  | -0.9426522510 | -0.4670861773 |
| 103 | -0.8090169944 | 0.9256147931  | -0.1261556613 | -1.0444364440 |
| 104 | -0.8090169944 | -0.9256147931 | 0.1261556613  | 1.0444364440  |
| 105 | -0.8090169944 | 1.1441228050  | 0.6605596176  | -0.4670861773 |
| 106 | -0.8090169944 | 1.1441228050  | -0.6605596176 | 0.4670861773  |
| 107 | -0.8090169944 | 1.2791681840  | 0.3302798088  | 0.4670861773  |
| 108 | 1.3090169944  | 0.3535533906  | -0.8646837645 | 0.1784110442  |
| 109 | 1.3090169944  | 0.5720614029  | -0.0779684865 | 0.7557613103  |
| 110 | -1.3090169944 | 0.2185080123  | -0.5344039563 | -0.7557613103 |
| 111 | -1.3090169944 | -0.2185080123 | 0.5344039563  | 0.7557613103  |
| 112 | -1.3090169944 | 0.3535533906  | 0.4564354698  | -0.7557613103 |
| 113 | -1.3090169944 | 0.3535533906  | -0.8646837645 | 0.1784110442  |
| 114 | -1.3090169944 | -0.3535533906 | -0.4564354698 | 0.7557613103  |
| 115 | -1.3090169944 | 0.5720914029  | 0.7385281032  | 0.1784110442  |
| 116 | -1.3090169944 | 0.5720614029  | -0.0779684865 | 0.7557613103  |
| 117 | -1.3090169944 | -0.5720614029 | -0.7385281032 | -0.1784110442 |
| 118 | -1.3090169944 | 0.9256147931  | -0.1261556613 | -0.1784110442 |
| 119 | -1.3090169944 | -0.9256147931 | 0.1261556613  | 0.1784110442  |
| 120 | -1.6180339891 | 0.            | 0.            | 0.            |

# The Joint Relations of the Vertices of the 600-Cell

| | |
|---|---|
| 1(2,3,8,9,17,18,23,24,29,38,40,99) | 33(8,12,15,17,20,21,29,31,34,48,49,119) |
| 2(1,5,7,8,18,20,24,26,28,29,35,38) | 34(10,14,15,19,20,21,30,32,33,49,117,119) |
| 3(1,6,9,17,23,37,40,64,79,85,96,99) | 35(2,5,7,22,24,26,38,41,44,51,54,56) |
| 4(5,10,16,19,22,25,30,36,41,63,66,97) | 36(4,22,25,39,41,51,55,63,66,77,91,108) |
| 5(2,4,7,10,16,20,22,26,28,30,35,41) | 37(3,6,23,40,42,50,53,57,64,76,79,90) |
| 6(3,11,12,17,23,31,37,42,64,76,82,96) | 38(1,2,8,23,24,26,35,40,46,53,54,56) |
| 7(2,5,16,18,22,24,28,35,44,71,73,103) | 39(25,36,45,52,55,61,63,69,77,78,89,108) |
| 8(1,2,12,15,17,20,23,26,29,33,38,46) | 40(1,3,9,23,24,37,38,53,54,57,79,92) |
| 9(1,3,18,24,40,73,74,79,85,92,99,105) | 41(4,5,10,22,25,26,35,36,47,51,55,56) |
| 10(4,5,14,15,19,20,25,26,30,34,41,47) | 42(6,11,12,23,27,37,43,46,50,53,58,76) |
| 11(6,12,13,21,27,31,42,43,67,76,82,104) | 43(11,13,27,42,45,50,52,58,67,76,78,88) |
| 12(6,8,11,15,17,21,23,27,31,33,42,46) | 44(7,22,24,35,51,54,59,71,73,81,91,92) |
| 13(11,14,21,27,32,43,45,67,69,78,84,104) | 45(13,14,25,27,39,43,47,52,55,58,69,78) |
| 14(10,13,15,19,21,25,27,32,34,45,47,69) | 46(8,12,15,23,26,27,38,42,47,53,56,58) |
| 15(8,10,12,14,20,21,26,27,33,34,46,47) | 47(10,14,15,25,26,27,41,45,46,55,56,58) |
| 16(4,5,7,22,28,30,66,71,97,102,103,110) | 48(17,29,31,33,49,96,99,111,112,115,119,120) |
| 17(1,3,6,8,12,23,29,31,33,48,96,99) | 49(20,28,29,30,33,34,48,110,112,117,119,120) |
| 18(1,2,7,9,24,28,29,73,99,103,105,112) | 50(37,42,43,52,53,57,58,60,76,88,90,109) |
| 19(4,10,14,25,30,32,34,63,69,97,100,117) | 51(22,35,36,41,44,54,55,56,59,60,91,108) |
| 20(2,5,8,10,15,26,28,29,30,33,34,49) | 52(39,43,45,50,55,58,60,78,88,89,108,109) |
| 21(11,12,13,14,15,27,31,32,33,34,104,119) | 53(23,37,38,40,42,46,50,54,56,57,58,60) |
| 22(4,5,7,16,35,36,41,44,51,66,71,91) | 54(24,35,38,40,44,51,53,56,57,59,60,92) |
| 23(1,3,6,8,12,17,37,38,40,42,46,53) | 55(25,36,39,41,45,47,51,52,56,58,60,108) |
| 24(1,2,7,9,18,35,38,40,44,54,73,92) | 56(26,35,38,41,46,47,51,53,54,55,58,60) |
| 25(4,10,14,19,36,39,41,45,47,55,63,69) | 57(37,40,50,53,54,59,60,79,90,92,94,109) |
| 26(2,5,8,10,15,20,35,38,41,46,47,56) | 58(27,42,43,45,46,47,50,52,53,55,56,60) |
| 27(11,12,13,14,15,21,42,43,45,46,47,58) | 59(44,51,54,57,60,81,91,92,93,94,108,109) |
| 28(2,5,7,16,18,20,29,30,49,103,110,112) | 60(50,51,52,53,54,55,56,57,58,59,108,109) |
| 29(1,2,8,17,18,20,28,33,48,49,99,112) | 61(39,62,63,68,69,77,78,83,84,89,98,100) |
| 30(4,5,10,16,19,20,28,34,49,97,110,117) | 62(61,65,67,68,78,80,84,86,88,89,95,98) |
| 31(6,11,12,17,21,33,48,82,96,104,111,119) | 63(4,19,25,36,39,61,66,69,77,83,97,100) |
| 32(13,14,19,21,34,69,84,100,104,114,117,119) | 64(3,6,37,65,70,76,79,82,85,90,96,101) |
| 65(62,64,67,70,76,80,82,86,88,90,95,101) | 93(59,68,72,75,77,80,81,89,91,94,108,109) |
| 66(4,16,22,36,63,71,72,77,83,91,97,102) | 94(57,59,70,74,75,79,80,81,90,92,93,109) |
| 67(11,13,43,62,65,76,78,82,84,88,95,104) | 95(62,65,67,82,84,86,98,101,104,111,114,116) |
| 68(61,62,72,75,77,80,83,86,89,93,98,106) | 96(3,6,17,31,48,64,82,85,99,101,111,115) |
| 69(13,14,19,25,32,39,45,61,63,78,84,100) | 97(4,16,19,30,63,66,83,100,102,110,113,117) |
| 70(64,65,74,75,79,80,85,86,90,94,101,107) | 98(61,62,68,83,84,86,95,100,106,113,114,116) |
| 71(7,16,22,44,66,72,73,81,87,91,102,103) | 99(1,3,9,17,18,29,48,85,96,105,112,115) |
| 72(66,68,71,75,77,81,83,87,91,93,102,106) | 100(19,32,61,63,69,83,84,97,98,113,114,117) |
| 73(7,9,18,24,44,71,74,81,87,92,103,105) | 101(64,65,70,82,85,86,95,96,107,111,115,116) |
| 74(9,70,73,75,79,81,85,87,92,94,105,107) | 102(16,66,71,72,83,87,97,103,106,110,113,118) |
| 75(68,70,72,74,80,81,86,87,93,94,106,107) | 103(7,16,18,28,71,73,87,102,105,110,112,118) |
| 76(6,11,37,42,43,50,64,65,67,82,88,90) | 104(11,13,21,31,32,67,82,84,95,111,114,119) |
| 77(36,39,61,63,66,68,72,83,89,91,93,108) | 105(9,18,73,74,85,87,99,103,107,112,115,118) |
| 78(13,39,43,45,52,61,62,67,69,84,88,89) | 106(68,72,75,83,86,87,98,102,107,113,116,118) |
| 79(3,9,37,40,57,64,70,74,85,90,92,94) | 107(70,74,75,85,86,87,101,105,106,115,116,118) |
| 80(62,65,68,70,75,86,88,89,90,93,94,109) | 108(36,39,51,52,55,59,60,77,89,91,93,109) |
| 81(44,59,71,72,73,74,75,87,91,92,93,94) | 109(50,52,57,59,60,80,88,89,90,93,94,108) |
| 82(6,11,31,64,65,67,76,95,96,101,104,111) | 110(16,28,30,49,97,102,103,112,113,117,118,120) |
| 83(61,63,66,68,72,77,97,98,100,102,106,113) | 111(31,48,82,95,96,101,104,114,115,116,119,120) |
| 84(13,32,61,62,67,69,78,95,98,100,104,114) | 112(18,28,29,48,49,99,103,105,110,115,118,120) |
| 85(3,9,64,70,74,79,96,99,101,105,107,115) | 113(83,97,98,100,102,106,110,114,116,117,118,120) |
| 86(62,65,68,70,75,80,95,98,101,106,107,116) | 114(32,84,95,98,100,104,111,113,116,117,119,120) |
| 87(71,72,73,74,75,81,102,103,105,106,107,118) | 115(48,85,96,99,101,105,107,111,112,116,118,120) |
| 88(43,50,52,62,65,67,76,78,80,89,90,109) | 116(86,95,98,101,106,107,111,113,114,115,118,120) |
| 89(39,52,61,62,68,77,78,80,88,93,108,109) | 117(19,30,32,34,49,97,100,110,113,114,119,120) |
| 90(37,50,57,64,65,70,76,79,80,88,94,109) | 118(87,102,103,105,106,107,110,112,113,115,116,120) |
| 91(22,36,44,51,59,66,71,72,77,81,93,108) | 119(21,31,32,33,34,48,49,104,111,114,117,120,) |
| 92(9,24,40,44,54,57,59,73,74,79,81,94) | 120(48,49,110,111,112,113,114,115,116,117,118,119) |